\documentclass[preprint,12pt]{elsarticle}
\usepackage{makeidx}
\usepackage{lineno,hyperref}
\usepackage{graphicx}
\usepackage{subfigure}
\usepackage{amsmath}
\usepackage{amsfonts}
\usepackage{captcont}
\modulolinenumbers[5]

\newtheorem{thm}{Theorem}
\newtheorem{conj}{Conjecture}

\newtheorem{lemma}{Lemma}

\usepackage{amssymb}

\journal{Journal Name}

\begin{document}

\begin{frontmatter}



\title{Development of a New Spectral Collocation Method Using Laplacian Eigenbasis for Elliptic Partial Differential Equations in an Extended Domain}


\author[mymainaddress]{P. Y. Wu}
\author[mymainaddress]{C. R. Kao}


\author[mysecondaryaddress,add3]{Tony W.H. Sheu\corref{mycorrespondingauthor}}
\cortext[mycorrespondingauthor]{Tony W.H. Sheu}
\ead{twhsheu@ntu.edu.tw}

\address[mymainaddress]{Department of Material Science and Engineering, National Taiwan University}
\address[mysecondaryaddress]{Department of Engineering Science and Ocean Engineering, National Taiwan University}
\address[add3]{Institute of Applied and Mathematical Science, National Taiwan University}

\begin{abstract}
The recent development of spectral method has been praised for its high-order convergence in simulating complex physical problems. The combination of embedded boundary method and spectral method becomes a mainstream way to tackle geometrically complicated problems. However, the convergence is deteriorated when embedded boundary strategies are employed. Owing to the loss of regularity, in this paper we propose a new spectral collocation method which retains the regularity of solutions to solve differential equations in the case of complex geometries. The idea is rooted in the basis functions defined in an extended domain, which leads to a useful upper bound of the Lebesgue constant with respect to the Fourier best approximation. In particular, how the stretching of the domain defining basis functions affects the convergence rate directly is detailed. Error estimates chosen in our proposed method show that the exponential decay convergence for problems with analytical solutions can be retained.  Moreover, two-dimensional Poisson equations and convection-diffusion equations with simple and complex geometrical domains will be simulated. The predicted results justify the advantages of applying our method to tackle geometrically complicated problems. 
\end{abstract}


\end{frontmatter}


\section{Introduction}

Spectral method has been well-known for its high accuracy in solving ordinary or partial differential equations (ODE or PDEs) in simple domains \cite{HAIDVOGEL1979167,CANUTO2007,Boyd2001}. A differential equation subject to periodic boundary conditions can be solved by spectral method to yield exponential decay convergence if the solution is analytic \cite{Gottlieb1977}. Recently, in the framework of spectral methods, boundary conditions of different kinds can be tackled by several proposed methods which can render high accuracy \cite{LI2009,Orszag1979273,PATERA1984468}. In view of the development of modern numerical schemes, spectral method is quite effective for dealing with physically complicated but geometrically simple problem. However, in practical applications, the predicted accuracy can be greatly deteriorated in geometrical complex domain. To retain high accuracy is still a crucial issue in the application of spectral method.

One promising way of tackling complex physical domain is the embedded boundary (EB) method \cite{Zhilin2006,FEDKIW1999457,Angot1999,STEIN2016252,LUI2009541}. EB method enables calculation of PDE solutions in geometrically complicated domains by simply extending computational domains which are suitable for discretization in structured grids. Most of the methods belonging to the category of EB method, such as Fourier spectral methods \cite{SABETGHADAM200955,CHEN1998147}, enjoy the nice feature of the high-order space discretization. Unfortunately, solutions of high continuous differentiability for differential equations may lose the smoothness when they are extended to the entire computational domain \cite{Beyer1992,Vega2007,STEIN2016252}. The degenerated regularity may lead to low-order convergence in the application of spectral method. Given that the high accuracy nature of spectral method is destroyed by the deteriorated regularity, many constructive methods have been proposed, including Immersed Interface \cite{Zhilin2006}, Fourier Continuation (FC) \cite{ALBIN20116248,BRUNO20102009,LYON20103358}, Active Penalty (AP) \cite{Shirokoff2015}, Immersed Boundary Smooth Extension (IBSE) methods \cite{STEIN2016252}, etc. The basic ideas of them are mainly rooted in grid refinement or smooth extension of functions near the boundary of physical domains. However, the convergence rates of these spectral methods are still restricted by their built in strategies. For example, the convergence rate of IBSE method is limited by the smoothness of $\delta$-function \cite{STEIN2016252}.

In this paper, a collocation method is proposed to resolve the difficulties of getting high accuracy when the complex physical domain is encountered. To the best of our knowledge, nearly all the spectral collocation methods can be categorized into the following groups in view of the applied basis functions: (i) eigenfunctions generated by solving the singular Sturm-Liouville problems and (ii) Fourier series.
In both basis functions, the accuracy of the scheme depends strongly on the collocation points. Flexibility of grid points to fit the boundary of geometrically complicated domains is therefore highly restricted, thereby motivating us to develop a new method through properly chosen collocation points. Our goal is that the convergence rate is less deteriorated. 

Different from other collocation methods, the basis generated by solving Laplace eigenvalue problem, subject to homogenous boundary conditions (say $u=0$ on the boundary), in a rectangular domain $R$ containing the computational domain $\Omega$ will be employed. We denote the process of constructing a basis in a larger domain containing the computational domain by \textit{stretching}. One can easily find that the basis functions are constituted of a sequence of sine functions, which are orthogonal and complete in $L^{2}(R)$. The necessity of stretching basis functions can be explained quantitatively as follows. Firstly, to deal with several different kinds of boundary conditions, the stretching of the basis functions is necessary since the basis functions are all zero on the boundary of $R$. Secondly, the size of $R$ can directly influence the convergence rate. The choice of $R$ is the core issue in our study and it will be discussed in Section 3. The basis constructed by our proposed stretching strategy yields an excellent convergence behavior as the traditional spectral collocation methods applied to problems defined in simple domains if the rectangular domain R is chosen properly. Inspired by the theory of trigonometric interpolation, we are led to find that the convergence with exponential decay for analytic solution of PDE can be retained by performing a  small perturbation to the uniform collocation points. Therefore, we can simply relocate collocation points to approximate the geometrically complex problem without much loss of the accuracy. To sum up, the new method we propose is applicable to simple geometry cases. Moreover, this method performs much better than most of the methods proposed in the past for solving geometrically complex problems.

The paper is organized as follows. The numerical model we propose will be given in Section 2. In Section 3, the error estimates for one-dimensional cases are given. It is found that that the accuracy is closely related to the stretching process. Numerical implementations of the proposed method to one- and two-dimensional problems with simple physical domains will be performed in Section 4. Numerical results for solving two-dimensional problems with complicated domain are presented in Section 5. We will then reiterate the new contribution of our proposed method and briefly describe future work.

\section{Methods}

We firstly consider a $d-$dimensional rectangular domain $R_{0} = (0,L_{1})\times(0,L_{2})\times\cdots\times (0,L_{d})$ and a larger rectangular domain $R=(-\delta_{1},L_{1}+\delta_{1})\times (-\delta_{2},L_{2}+\delta_{2})\times\cdots\times (-\delta_{d},L_{d}+\delta_{d})$, where $L_{i}$ and $\delta_{i}$, $i=1,\ldots ,d$ are some constants. The modeling problem for the description of our numerical method is expressed as
\begin{subequations}\label{eq:1}
	\begin{align}
	-\Delta u = f \ \ \ \ 
	&\text{in} \ R_{0}, \\
	u = g \ \ \ \ 
	&\text{on} \ \partial R_{0}.
	\end{align}
\end{subequations}
Given that $f$ is smooth and $g$ is smooth as well in the sense that there exists an extension $\tilde{g}$ of $g$ such that $\tilde{g}$ is smooth in $R_{0}$ and $\tilde{g}=g$ on $\partial R_{0}$.
Basis functions for the approximation of Eq. (\ref{eq:1}) are defined by Laplacian eigenvalue problem subject to homogeneous boundary conditions in $R$ and the restrictions in $R_{0}$ are employed. That is, basis functions are generated by solving the equation
\begin{subequations}\label{eq:2}
	\begin{align}
	-\Delta w = \lambda w \ \ \ \ 
           &\text{in} \ R, \\
	w = 0 \ \ \ \ 
           &\text{on} \ \partial R.
	\end{align}
\end{subequations}
The eigenfunctions of Eq. (\ref{eq:2}) can be written in the form of
\begin{equation}\label{eq:3}
	w_{j_{1}j_{2}\ldots j_{d}}(x_{1},x_{2},\ldots ,x_{d})=\frac{2^{\frac{d}{2}}}{\sqrt{\prod_{k=1}^{d}(L_{k}+2\delta_{k})}}\prod_{k=1}^{d}\sin \left( \frac{j_{k}\pi (x_{k}+\delta_{k})}{L_{k}+2\delta_{k}}\right)
\end{equation}
associated with the eigenvalues given below
\begin{equation}\label{eq:4}
	\lambda_{j_{1}j_{2}\ldots j_{d}}=\sum_{k=1}^{d}\frac{j_{k}^{2}\pi^{2}}{(L_{k}+2\delta_{k})^{2}} ,
\end{equation}
where $j_{1},j_{2},\ldots ,j_{d}\in \mathbb{N}$. Thanks to the construction of basis functions, Eq. (\ref{eq:1}) can be converted to the problem: Find $a_{i}$ such that
\begin{subequations}\label{eq:5}
	\begin{align}
	\sum_{i=1}^{\infty}a_{i}\lambda_{i}w_{i}=f \ \ \ \ 
	&\text{in} \ R_{0}, \\
	\sum_{i=1}^{\infty}a_{i}w_{i}=g \ \ \ \ 
	&\text{on} \ \partial R_{0}.
	\end{align}
\end{subequations}
Here the eigenbasis is numbered by the index $i$. The eigenvalues $\lambda_{j}$ are subject to the constraint given by
\begin{equation}\label{eq:6}
	\lambda_{1}<\lambda_{2}\leq \lambda_{3} \leq \cdots ,
\end{equation}
where $\lambda_{i}\rightarrow \infty$ as $i\rightarrow \infty$.

The collection of all eigenfunctions in Eq. (\ref{eq:2}), which is
\begin{equation*}
	W := \lbrace w_{j_{1}j_{2}\ldots j_{d}}\rbrace_{j_{1},j_{2},\ldots ,j_{d}\in\mathbb{N}},
\end{equation*}
is an orthonormal basis in $L^{2}(R)$. Notice that the restriction on $R_{0}$ of $W$, say $W_{R_{0}}$, spans $L^{2}(R_{0})$. However, $W_{R_{0}}$ is not linearly independent since, for example, the function vanishing in $\overline{R_{0}}$ but taking an arbitrary nonzero constant in $R-\overline{R_{0}}$ is $L^{2}$-integrable. Fortunately, the finite collection of $w_{j_{1}j_{2}\ldots j_{d}}|_{R_{0}}$ is linearly independent. The availability of the numerical implementation of the basis $W_{R_{0}}$ is therefore ensured. Linear independence of the finite collection from the elements in $W_{R_{0}}$ will be proved in Section 3.

\subsection{Space discretizations in rectangular domains}

For the space discretization of Eq. (\ref{eq:1}), we firstly consider an approximated solution of the form
\begin{equation}\label{eq:7}
	u_{n}=\sum_{j_{1},j_{2},\ldots ,j_{d}=1}^{n}a_{j_{1}j_{2}\ldots j_{d}}w_{j_{1}j_{2}\ldots j_{d}},
\end{equation}
where $a_{j_{1}j_{2}\ldots j_{d}}$'s are some real coefficients. Our goal is to find these unknown coefficeints $a_{j_{1}j_{2}\ldots j_{n}}$ so that $\|u-u_{n}\|_{L^{2}(R_{0})}$ is minimized. To reach this goal, the uniform grid which has $n$ grid points in each direction of the rectangular $C$ is chosen. We denote the collection of uniform grid points by X, where
\begin{equation}\label{eq:8}
\begin{aligned}
	&X=\Big\{ x=(x^{1}_{j_{1}},x^{2}_{j_{2}},\ldots ,x^{d}_{j_{d}}) \big| \ x^{k}_{1}=0  , \  x^{k}_{n}=L_{k}, \\
  &|x^{k}_{j_{k}+1}-x^{k}_{j_{k}}| = \frac{L_{k}}{n-1}  ,\ j_{1},j_{2},\ldots j_{d} = 1,\ldots , n, \ k=1,\ldots , d \Big\}.
\end{aligned}
\end{equation}
By an abuse of notation, we denote the uniform grid points by $ X = \lbrace x_{i} \rbrace $ and the finite collection of basis functions by $W_{R_{0}}=\lbrace w_{j}\rbrace $. The smallest rectangle formed by $2^{d}$ points in the uniform grid system is called an "element". Under the circumstances, the cardinality of $X$ and $W_{R_{0}}$ are identical. We may group grid points into basis functions to get a system of $n^{d}$ linear equations
\begin{subequations}\label{eq:9}
	\begin{align}
	\sum_{j=1}^{n^{d}}a_{j}\lambda_{j}w_{j}(x_{i})=f(x_{i}) \ \ \ \ 
	&\text{for}\ x_{i}\in R_{0}, \\
	\sum_{j=1}^{n^{d}}a_{j}w_{j}(x_{i})=g(x_{i}) \ \ \ \ 
	&\text{for}\ x_{i}\in \partial R_{0},
	\end{align}
\end{subequations}
where $i=1,\ldots , n^{d}$. We can obtain an approximated solution to Eq. (\ref{eq:1}) by calculating $a_{i}$ from Eq. (\ref{eq:9}).

\subsection{Space discretizations in irregular domains}

\begin{figure}
	\centering
	\includegraphics[width=0.6\linewidth]{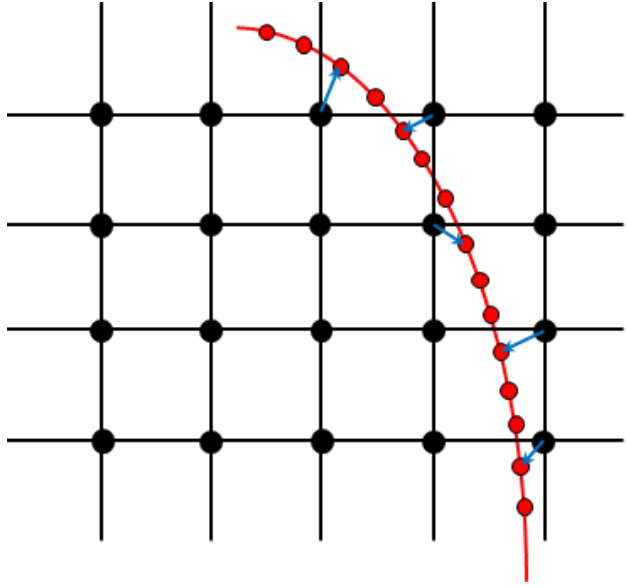}
	\caption{Schematic of the way how to move the grid points on the boundary.}
	\label{schematic_moved}
\end{figure}

To explain the idea of approximating PDEs in an irregular domain, we restrict ourselves to two-dimensional cases. Let $\Omega\subset \mathbb{R}^{2}$ be a smooth domain. Without loss of generality, we may consider that $\Omega$ is contained in a rectangular domain $R_{0}=(0,L_{1})\times (0,L_{2})$. Basis functions and their associated eigenvalues are defined by Eq. (\ref{eq:3}) and (\ref{eq:4}). Similar to Eq. (\ref{eq:1}), the problem under investigation in irregular domain $\Omega$ is given by
\begin{subequations}\label{eq:10}
	\begin{align}
	-\Delta u = f \ \ \ \ 
	&\text{in} \ \Omega, \\
	u = g \ \ \ \ 
	&\text{on} \ \partial \Omega.
	\end{align}
\end{subequations}
Let $X$ be a two-dimensional uniform grid of $R_{0}$ defined by Eq. (\ref{eq:8}), and $\Gamma$ a closed curve which defines $\partial \Omega$. The uniform grid points are not compatible with $\partial \Omega$ in general. Therefore, we move the grid points near the closed curve $\Gamma$ to some points exactly on the curve $\Gamma$. In the numerical implementation, we may define $\Gamma$ as the following:
\begin{equation}\label{eq:11}
	\Gamma := \lbrace \gamma (t)=(x(t),y(t)) \big| t\in [0,1], \ x(0)=x(1),\ y(0)=y(1) \rbrace .
\end{equation}
In addition, the winding number of $\Gamma$ is $1$. The grid points of $\Gamma$, denoted by $X_{\Gamma}$, can be expressed as
\begin{equation}\label{eq:12}
	X_{\Gamma} := \lbrace \gamma(t_{i})=(x(t_{i}),y(t_{i})) \big| t_{i}=\frac{i}{n_{c}}, \ i = 1,\ldots ,n_{c} \rbrace ,
\end{equation}
where $n_{c}$ is the number of grid points on $\Gamma$. For a pair of grid points $\lbrace\gamma(t_{i}), \gamma(t_{i+1})\rbrace$ on $\Gamma$ passing through an edge of an element, the point of $X$ on the edge which is near $\Gamma$ will be moved to $\gamma(t_{i})$, see Fig. \ref{schematic_moved} and \ref{35x35}. Furthermore, The boundary condition shown in Eq. (\ref{eq:10}b) will be imposed on those points being relocated along $\partial \Omega$. The system of linear equations for the approximated solution to Eq. (\ref{eq:10}) is given by
\begin{subequations}\label{eq:13}
	\begin{align}
	\sum_{j=1}^{n^{2}}a_{j}\lambda_{j}w_{j}(x_{i})=f(x_{i}) \ \ \ \ 
	&\text{for}\ x_{i}\in R_{0}, \\
	\sum_{j=1}^{n^{2}}a_{j}w_{j}(x_{i})=g(x_{i}) \ \ \ \ 
	&\text{for}\ x_{i} \ \text{being moved}.
	\end{align}
\end{subequations}

\begin{figure}
	\centering
    \subfigure[The $35\times 35$ uniform grid.]{
    \label{35x35}
    \includegraphics[width=0.62\linewidth]{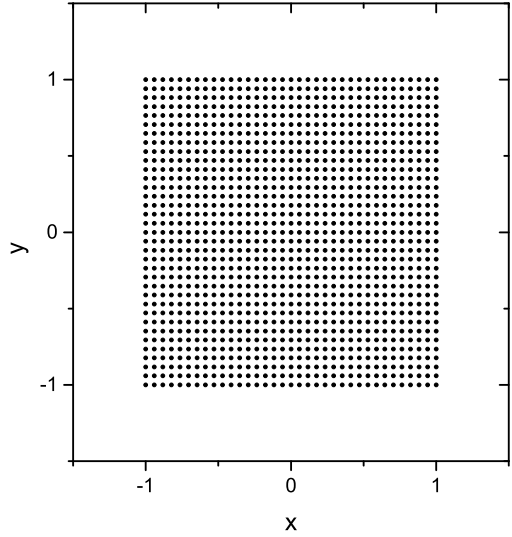}
    }
    \subfigure[The grid being relocated.]{
    \label{35x35_moved}
    \includegraphics[width=0.62\linewidth]{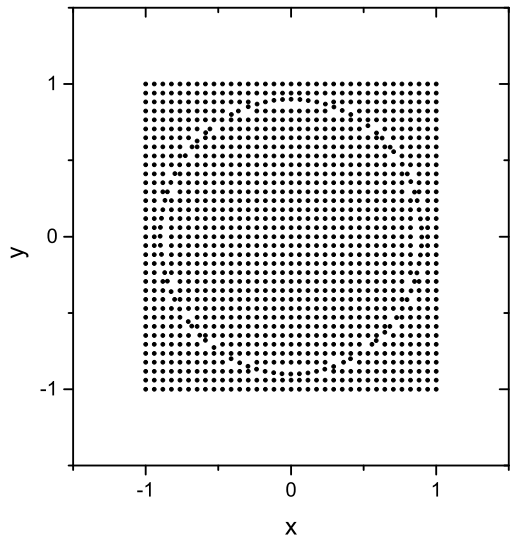}
    }
    \caption{An example chosen to show how the grid points are moved. In this case, we move the grid points near the boundary of the disk $B_{0.9}((0,0))$ to the boundary of the disk.}
    \label{35x35}
\end{figure}

\subsection{Generalization}
In many cases of elliptic problem, first order differential terms are involved in PDEs. For example, convection-diffusion-reaction equations are frequently applied to model transport phenomena in the real world. Refering to Eq. (\ref{eq:7}), the first order partial differential of $u_{n}$ can be written as
\begin{equation}\label{eq:14}
	\frac{\partial u_{n}}{\partial x_{i}}=\sum_{j_{1},j_{2},\ldots ,j_{d}=1}^{n}a_{j_{1}j_{2}\ldots j_{d}}\frac{\partial}{\partial x_{i}}w_{j_{1}j_{2}\ldots j_{d}}.
\end{equation}
It is noted that
\begin{equation}\label{eq:15}
	\frac{\partial}{\partial x_{i}}w_{j_{1}j_{2}\ldots j_{d}}=\frac{2^{\frac{d}{2}}j_{i}\pi \cos\left( \frac{j_{i}\pi (x_{i}+\delta_{i})}{L_{i}+2\delta_{i}}\right)}{(L_{i}+2\delta_{i})\sqrt{\prod_{k=1 }^{d}(L_{k}+2\delta_{k})}} \prod_{k=1, k\neq i}^{d}\sin \left( \frac{j_{k}\pi (x_{k}+\delta_{k})}{L_{k}+2\delta_{k}}\right).
\end{equation}
The collection of $\frac{\partial}{\partial x_{i}}w_{j_{1}j_{2}\ldots j_{d}}$ forms a linearly independent set, which is also linearly independent of $W$. Therefore, we can simply add some terms of the form Eq. (\ref{eq:14}) to the algebraic equations Eq. (\ref{eq:9}) or (\ref{eq:13}) for solving the approximated solution to the PDE with the first order differential terms. 

\section{Convergence Estimates}
In this section, we are going to discuss the convergence of the approximated solution and the choice of extended domain $R$ that shall affect the convergence rate of our proposed method. The modeling problems discussed here are the convergence estimates on Poisson's equation and convection-diffusion equations. For the sake of simplicity, the one-dimensional case is considered. Moreover, all the arguments about the interpolation are under the assumption of the uniform grid $X := \lbrace \frac{(j-1)L}{N-1} \rbrace_{j=1,\ldots ,N}$, where $N$ is a integer larger than $3$ and $L$ is a positive number. All the functions mentioned in this paper are assumed to be bounded everywhere.
\subsection{One-dimensional Poisson's equation}
Let $I=(0,L)$, where $L$ is a positive constant. The modeling problem is given by
\begin{subequations}\label{eq:16}
	\begin{align}
		-u_{xx}=f \ \ \ \ \text{in}\ I \\
	   	u(0)=a, \  u(L)=b,
	\end{align}
\end{subequations}
where $f\in C^{m}(I)$.
Let $\delta$ be a positive number, $w_{j}$ the basis functions defined in Eq. (\ref{eq:3}) in $\widetilde{I}:=(-\delta , L+\delta )$, $u_{N}$ the approximated solution determined by Eq. (\ref{eq:9}), where $N$ is a positive integer. Our goal is to find the relation between the error $\| u-u_{N} \|_{L^{2}(I)}$ or $\| u-u_{N} \|_{L^{\infty}(I)}$ and $N$.

Before performing the convergence estimation, we shall show that the interpolation process determined by the assembled matrix in Eq. (\ref{eq:9}) is well-defined. To be specific, we need to check whether the system given in Eq. (\ref{eq:9}) is nonsingular.

\begin{lemma}
Let $w_{j}$'s be the basis functions generated from Eq. (\ref{eq:2}) for the dimension $d=1$. The finite collection of $w_{j}$'s restricted on $I=(0,L)$ is linearly independent. \\
\textbf{Proof.} Assume that $\sum_{j=1}^{N}a_{j}w_{j}|_{I}=0$, where $a_{j}$'s are some real numbers, $N$ is a positive integer. We claim that $\sum_{j=1}^{N}a_{j}w_{j}=0$ in $\widetilde{I}:= (-\delta , L+\delta )$, where $\delta$ is a positve number. Otherwise, there exists $x_{0}$ in $\widetilde{I}-I$ such that $\sum_{j=1}^{N}a_{j}w_{j}(x_{0})\neq 0$. Let $p(x)=\sum_{j=1}^{N}a_{j}w_{j}(x)$. Since $p$ is analytic in $\widetilde{I}$, we have $f(x_{0})=\sum_{j=0}^{\infty}\frac{p^{(j)}(x)}{j!}(x_{0}-x)^{j}$ for all $x$ in $I$. However, $p^{(j)}(x)=0$ for all $x$ in I, we have $p(x_{0})=0$, thereby leading to a contradiction. Since $\lbrace w_{j}\rbrace_{j=1}^{\infty}$ is linearly independent, we have $a_{j}=0$ for $j=1,\ldots , N$. This completes the proof.
\end{lemma}

Now, we have to show that the error estimates for one-dimensional Poisson's equation can be made in the interval $\widetilde{I}$. Let $\widetilde{f}$ be the $C^{m}$ extension of $f$ such that $\widetilde{f}=f$ in $I$ and $\widetilde{f}$ has a compact support in $\widetilde{I}-I$. Also, let $\widetilde{u}$ be the extension of the solution $u$ for Eq. (\ref{eq:16}) defined by
\begin{subequations}\label{eq:17}
	\begin{align}
		-u_{xx}=\widetilde{f}\ \ \ \  \text{in} \ \widetilde{I}, \\
		\widetilde{u}(-\delta )=\widetilde{u}(L+\delta )=0, \ \widetilde{u}(0)=a, \ \widetilde{u}(L)=b.
	\end{align}
\end{subequations}
The existence of $\widetilde{u}$ and the associated $\widetilde{f}$ can be ensured from the classical extension theorem.

By applying Poincar\'{e}'s inequality, there exists a constant $C_{P}$ that is independent of $N$ such that
\begin{equation}\label{eq:18}
	| u-u_{N} |_{L^{2}(I)} \leq |\widetilde{u}-\widetilde{u}_{N} |_{L^{2}(\widetilde{I})} \leq  C_{P}|\partial_{xx}(\widetilde{u}-\widetilde{u}_{N}) |_{L^{2}(\widetilde{I})} = C_{P}| \widetilde{f}-f_{N} |_{L^{2}(\tilde{I})},
\end{equation} 
where $f_{N}$ is the trigonometric interpolation by $\lbrace w_{j}\rbrace_{j=1}^{N}$. Note that $f_{N}(x_{i})=f(x_{i})$ for $x_{i}$ described in Eq. (\ref{eq:9}). So far, our problem is simplified to the estimation of $\|\widetilde{f}-f_{N}\|_{L^{2}(\widetilde{I})}$. 
Since $\widetilde{f}\in L^{2}(\widetilde{I})$, $\widetilde{f}$ can be expressed as $\widetilde{f}=\sum_{j=1}^{\infty}b_{j}w_{j}$, where $b_{j}= (\widetilde{f},w_{j})_{L^{2}(\widetilde{I})}$. For the space $V_{N}:=span\lbrace w_{1},\ldots , w_{N}\rbrace$, we can see that $\sum_{j=1}^{N}b_{j}w_{j}$ is the best approximation for $\widetilde{f}$ in $V_{N}$. We shall observe the convergence rate of the best approximation with respect to $\widetilde{f}$.

\begin{lemma}
	Let $m$ be a positive integer. If $\widetilde{f}\in C^{m}_{0}(\widetilde{I})$ piecewisely, and all the derivatives of $\widetilde{f}$ vanish at the boundary of $\widetilde{I}$, then the $L^{2}$-error $\| \widetilde{f}-\sum_{j=1}^{N-1}b_{j}w_{j} \|_{L^{2}(\widetilde{I})}$ is $O(N^{-m-1})$. \\
\textbf{Proof.} Applying integration by parts, $b_{N}$ can be calculated as
\begin{equation*}
	\begin{aligned}
		b_{N} =& \frac{\sqrt{2(L+2\delta )}}{N\pi}\lbrace [ f(-\delta )-(-1)^{N}f(L+\delta )] \\
          	&+\frac{L+2\delta}{N\pi}\int_{-\delta}^{L+\delta}f^{'}(x)\cos\left(\frac{N\pi (x+\delta)}{L+2\delta}\right)dx \rbrace \\
		=& \frac{\sqrt{2(L+2\delta )}}{N\pi}\lbrace [ f(-\delta )-(-1)^{N}f(L+\delta )] \\
          	&-\frac{(L+2\delta )^{2}}{N^{2}\pi^{2}}\int_{-\delta}^{L+\delta}f^{''}(x)\sin\left(\frac{N\pi (x+\delta)}{L+2\delta}\right)dx \rbrace \\
		=& \cdots \\
		=& \sum_{\substack{k\geq 0, \\ k\ even}}^{m-1}\frac{\sqrt{2}(L+2\delta )^{\frac{1}{2}+k}}{(N\pi )^{k+1}}(-1)^{\frac{k}{2}}[f^{(k)}(-\delta )-(-1)^{N}f^{(k)}(L+\delta )] + O(N^{-m-1}).
	\end{aligned}
\end{equation*}
Since $\widetilde{f}$ and the derivatives of $f$ vanish at the boundary of $\widetilde{I}$, the proof is completed.
\end{lemma}
By Lemma 2, we find that the convergence rate of the best approximation behaves exponentially if $f$ is smooth. To justify our proposed method, we shall prove that the convergence rate of the approximation $f_{N}$ won't differ too much from the best approximation's.

Recalling Eq. (\ref{eq:9}), $f_{N}$ can be expressed in terms of $f(x_{i})$
\begin{equation}\label{eq:19}
	f_{N}(x)=\sum_{j=1}^{N}f(x_{j})l_{j}(x),
\end{equation}
where $l_{j}$'s are trigonometric Lagrange's interpolations. Let $z=x+\delta$, $z_{k}=x_{k}+\delta$ for $k=1,\ldots N$, $l_{j}$ can be expressed as \cite{Zygmund1959}
\begin{equation}\label{eq:20}
	l_{j}(z) = \frac{\sin\left( \frac{\pi z}{L+2\delta}\right) \prod_{k\neq j}\left( \cos\left(\frac{\pi z}{L+2\delta}\right)-\cos\left(\frac{\pi z_{k}}{L+2\delta}\right)\right)}{\sin\left( 	\frac{\pi z_{j}}{L+2\delta}\right) \prod_{k\neq j}\left( \cos\left(\frac{\pi z_{j}}{L+2\delta}\right)-\cos\left(\frac{\pi z_{k}}{L+2\delta}\right)\right)}.
\end{equation}
Introducing the Lebesgue function
\begin{equation}\label{eq:21}
	L(z)=\sum_{j=1}^{N}| l_{j}(z)|,
\end{equation}
the Lebesgue constant is
\begin{equation}\label{eq:22}
	\Lambda_{N}=\max_{z\in [0,L]}L(z).
\end{equation}
Now, the relation between $f_{N}$ and the best approximation can be written explicitly as
\begin{equation}\label{eq:23}
	| f-f_{N} |_{L^{2}(\widetilde{I})} \leq (1+ \Lambda_{N})| f-\sum_{j=1}^{N}b_{j}w_{j}|_{L^{2}(\widetilde{I})} .
\end{equation}
Therefore, the convergence rate of $f_{N}$ can be estimated by computing the values of $|l_{j}(x)|$ for each $j=1,\ldots , N$. 
Firstly, the denominator of $l_{j}(z)$ given by Eq. (\ref{eq:21}) can be estimated. Given $z_{j}$, $j=1,\ldots ,N$ as in Eq. (\ref{eq:21}), we have
\begin{equation}\label{eq:24}
	\left| \sin \left( \frac{\pi z_{j}}{L+2\delta}\right) \right| \geq \left| \sin \left( \frac{\pi \delta}{L+2\delta}\right) \right| ,
\end{equation}
and
\begin{equation}\label{eq:25}
	\begin{aligned}
		&\prod_{k\neq j}\left|\left( \cos\left(\frac{\pi z_{j}}{L+2\delta}\right)-\cos\left(\frac{\pi z_{k}}{L+2\delta}\right)\right) \right| \\
		&\geq \left| \cos\left(\frac{\pi \delta}{L+2\delta}\right)-\cos\left(\frac{\pi (\delta +h)}{L+2\delta}\right) \right|^{N-1} \\
		&= \left| 2\sin\left(\frac{\pi (\delta+\frac{1}{2}h)}{L+2\delta}\right)\sin\left(\frac{\pi (\delta-\frac{1}{2}h)}{L+2\delta}\right) \right|^{N-1},
	\end{aligned}
\end{equation}
where $h=\frac{L}{N-1}$. To find the upper bound of the remaining numerator terms of $l_{j}(z)$, we consider polynomials $\phi_{j}$, $j=1,\ldots ,N$, on a unit circle $S_{1}$ in the complex plane $\mathbb{C}$.
\begin{equation}\label{eq:26}
	\phi_{j}(\zeta )=\prod_{k\neq j} (\zeta - e^{\frac{i\pi z_{k}}{L+2\delta}}), \ \zeta\in S_{1}.
\end{equation}
Let $\zeta = e^{\frac{i\pi z}{L+2\delta}}$, where $z\in [-\delta , L+\delta ]$, we can get
\begin{equation}\label{eq:27}
	|\phi_{j}(\zeta )| \geq \prod_{k\neq j}\left|\left( \cos\left(\frac{\pi z}{L+2\delta}\right)-\cos\left(\frac{\pi z_{k}}{L+2\delta}\right)\right) \right|.
\end{equation}
Note that the right hand side of Eq. (\ref{eq:28}) is not smaller than the absolute value of the numerator of $l_{j}(z)$. The estimation can be continued by determining the upper bound of $|\phi_{j}(\zeta )|$. Upon calculating the distances between $\zeta$ to all roots of $\phi_{j}$, we have
\begin{equation}\label{eq:28}
	\begin{aligned}
		|\phi_{j}(\zeta )|\leq \left(2\cos \left(\frac{\delta\pi}{L+2\delta}\right)\right)\left(2\cos \left(\frac{2L\pi}{2(N-1)(L+2\delta)}\right)\right)\cdots \\
		\left(2\cos \left(\frac{3L\pi}{2(N-1)(L+2\delta)}\right)\right)\cdots \left(2\cos \left(\frac{(N-1)L\pi}{2(N-1)(L+2\delta)}\right)\right).
	\end{aligned}
\end{equation}
Therefore, we have
\begin{equation}\label{eq:29}
	|\phi_{j}(\zeta )|= 2^{N-1} \cos \left(\frac{\delta\pi}{L+2\delta}\right)\sin (2\omega\pi) \cdots \sin ((N-1)\omega\pi ),
\end{equation}
where $\omega = \frac{L}{2(N-1)(L+2\delta )}$. Combining Eq. (\ref{eq:25}), (\ref{eq:26}), (\ref{eq:28}) and (\ref{eq:30}), we have
\begin{equation}\label{eq:30}
	| l_{j}(z)| \leq 2^{N-1}\cot (\frac{\delta\pi}{L+2\delta}) \frac{\sin(2\omega\pi )\cdots \sin((N-1)\omega\pi )}{\left| 2\sin\left(\frac{\pi (\delta+\frac{1}{2}h)}{L+2\delta}\right)\sin\left(\frac{\pi (\delta-\frac{1}{2}h)}{L+2\delta}\right) \right|^{N-1}}.
\end{equation}
The above arguments can be integrated into a theorem:
\begin{thm}
	There exists a positive number $\delta$ used in Eq. (\ref{eq:31}) such that
\begin{equation}\label{eq:31}
	\Lambda_{N}\leq 2N\cot \left(\frac{\delta \pi}{L+2\delta} \right) ,
\end{equation}
where $\Lambda_{N}$ is the Lebesgue constant given in Eq. (\ref{eq:23}). 
\end{thm}

\subsection{One-dimensional convection-diffusion equation}
The following equation is then considered:
\begin{subequations}\label{eq:32}
	\begin{align}
	-u_{xx}+ku_{x}=f\ \ \ \text{in}\ I=(0,L), \\
    u(0)=a,\ u(L)=b.
    \end{align}
\end{subequations}
In the above, let $m\in \mathbb{N}\cup \lbrace 0\rbrace $, we assume that $k,f\in C^{m}(I)$. By the argument in Section 3.1, the extension of the solution to Eq. (\ref{eq:32}) can be also ensured:
\begin{lemma}
	Let $m\in \mathbb{N}\cup\lbrace 0\rbrace$, $u\in C^{m+2}(I)$, $k\in C^{m}(\widetilde{I})$, $f\in C^{m}(I)$ such that
		$$ -u_{xx}+ku_{x}=f \ \ \text{in} \ I, $$
		$$ u(0)=a, \ u(L)=b. $$
		Then, there exists $\widetilde{u}\in C^{m+2}_{0}(\widetilde{I})$, $\widetilde{f}\in C^{m}_{0}(\widetilde{I})$ satisfying
		$$ -\widetilde{u}_{xx}+k\widetilde{u}_{x}=\widetilde{f}\ \text{in} \ \ \widetilde{I}, $$
		$$ \widetilde{u}(-\delta)=\widetilde{u}(L+\delta)=0, $$
		where $\widetilde{f}=f$ in $(0,L)$ and $\widetilde{u}=u$ in $[0,L]$.
\end{lemma}
The estimation to the approximated solution of the convection-diffusion equations can be simplified to the estimation of the approximation of the forcing term. Hence, we have the following lemma:
\begin{lemma}
	Assuming that $u\in C_{0}^{2}(\widetilde{I})$, $k\in C(\widetilde{I})$, $f\in C_{0}(\widetilde{I})$, there is a constant $C=C(L,\delta , b)$ such that 
		$$|u-u_{N}|_{L^{2}(\widetilde{I})}\leq C|f-f_{N}|_{L^{2}(\widetilde{I})}, $$
		where
		$$f_{N}(x)=\sum_{i=1}^{N}\lambda_{i}a_{i}w_{i}(x)+ka_{i}w_{i,x}(x).$$
	\textbf{Proof.} Let $y=u-u_{N}$, $g=f-f_{N}$, we have
    \begin{equation}\label{eq:33}
	 	-y_{xx}+ky_{x}=g.
    \end{equation}
	Multiplying Eq. \ref{eq:33} by $-y_{x}$, we have
	$$ y_{xx}y_{x}-ky_{x}^{2}=-gy_{x}. $$
	Let $K(x) = \int_{-\delta}^{x}k(s)ds $, the following differential equation can be obtained:
	$$ \frac{d}{dx}(e^{-2K(x)}y_{x}^{2}) = -2gy_{x}e^{-2K(x)}. $$
	Therefore,
	$$  \frac{d}{dx}(e^{-2K(x)}y_{x}^{2}) \leq 2|g||y_{x}|e^{-2K(x)}, $$
    $$ e^{-2K(x)}y_{x}^{2} \leq \int_{-\delta}^{L+\delta} 2|g||y_{x}|e^{-2K(s)}ds \leq 2(L+2\delta )|g|_{L^{2}(\widetilde{I})} |y_{x}|_{L^{2}(\widetilde{I})}e^{2|K|_{L^{\infty}(\widetilde{I})}}.$$
	By integrating the terms on both sides of the above equation, we can obtain
	$$ |y_{x}|_{L^{2}(\widetilde{I})} \leq 2(L+2\delta )^{2}e^{4|K|_{L^{\infty}(\widetilde{I})}}|g|_{L^{2}(\widetilde{I})}. $$
	Since $|K(x)| \leq (L+2\delta )|k|_{L^{\infty}(\widetilde{I})}$ for all $x\in\widetilde{I}$, we have
	$$ |y_{x}|_{L^{2}(\widetilde{I})} \leq 2(L+2\delta )^{2} e^{4(L+2\delta)|k|_{L^{\infty}(\widetilde{I})}} |g|_{L^{2}(\widetilde{I})}. $$
    By applying Poincar\'{e}'s inequality, the proof can be completed.
\end{lemma}

By the theory of functional analysis, there is a positive number $\Lambda_{N}$ such that
\begin{equation}\label{eq:34}
	|f-f_{N}|_{L^{2}(\widetilde{I})} \leq (1+\Lambda_{N}) |f-f_{*}|_{L^{2}(\widetilde{I})},
\end{equation}
where $f_{*}$ is the best approximation in terms of the basis $\lambda_{j}w_{j}+kw_{j,x}$, $j=1,\ldots ,N$. However, the best approximation expressed in such basis is difficult to be obtained. We may find a function that bounds the forcing term and whose best approximation can be found more easily. Let $u_{N}$ be any kind of N-points interpolation of $u$ and $f_{N}$ be the associated approximated forcing term satisfying
\begin{equation}\label{eq:35}
	f_{N}=u_{N,xx}+ku_{N,x}\ \ \ \text{in}\ I.
\end{equation}
Combining Eq. (\ref{eq:35}) with Eq. (\ref{eq:32}), we have
\begin{equation*}
	\begin{split}
		|f-f_{N}|_{L^{2}(\widetilde{I})} & \leq |(u-u_{N})_{x}|_{L^{2}(\widetilde{I})} + |k|_{L^{\infty}(\widetilde{I})} |(u-u_{N})_{xx}|_{L^{2}(\widetilde{I})} \\
			& \leq ( C(L,\delta )+|k|_{L^{\infty}(\widetilde{I})}) |(u-u_{N})_{xx}|_{L^{2}(\widetilde{I})}.
	\end{split}
\end{equation*}
Therefore,
	$$ |f-f_{N}|_{L^{2}(\widetilde{I})} \leq (C(L,\delta )+|k|_{L^{\infty}(\widetilde{I})})(1+\Lambda_{N}) |(u-u_{N}^{*})_{xx}|_{L^{2}(\widetilde{I})}, $$
    where $u_{N}^{*}$ is the best approximation in terms of the basis $w_{j}$, $j=1,\ldots ,N$.
    
To express the best approximation $|(u-u_{N})_{xx}|_{L^{2}(\widetilde{I})}$ explicitly, we need the following lemma:
\begin{lemma}
	Let $h\in C_{0}^{m}(\widetilde{I})$, $h^{*}$ be the best approximation of $h$ in terms of the basis $w_{j}$, $j=1,\ldots ,N$, then we have
	\begin{equation}\label{eq:36}
    	|h-h^{*}|_{L^{2}(\widetilde{I})} \leq \frac{ (L+2\delta )^{m+2}}{N^{m+1}\pi^{m+1}} |h^{(m)}|_{L^{2}(\widetilde{I})}^{2},
	\end{equation}
    where $h^{(m)}$ is the $m$-th derivative of $h$. \\
    \textbf{Proof.} By Lemma 2, the $N-th$ Fourier coefficient $b_{N}$ of $h$ can be expressed either as
    $$ b_{N} = \pm \int_{-\delta}^{L+\delta}\frac{\sqrt{2}(L+2\delta )^{m+3/2}}{N^{m+1}\pi^{m+1}}h^{(m)}\cos\left( \frac{N\pi (x+\delta )}{L+2\delta}\right) dx, $$
	or
	$$ = \pm \int_{-\delta}^{L+\delta}\frac{\sqrt{2}(L+2\delta )^{m+3/2}}{N^{m+1}\pi^{m+1}}h^{(m)}\sin\left( \frac{N\pi (x+\delta )}{L+2\delta}\right) dx. $$
    Therefore, we have
    $$ |b_{N}|^{2} \leq \frac{ (L+2\delta )^{m+2}}{N^{m+1}\pi^{m+1}} |b_{N}^{(m)}|^{2}. $$
    Summing up $|b_{j}|^{2}$ for $j \geq N$, the following inequality can be derived
	$$ |\sum_{j=N}^{\infty}b_{j}w_{j}|_{L^{2}(\widetilde{I})}^{2}\leq \frac{ (L+2\delta )^{m+2}}{N^{m+1}\pi^{m+1}} \sum_{j=N}^{\infty} |b_{j}^{(m)}|^{2} \leq \frac{ (L+2\delta )^{m+2}}{N^{m+1}\pi^{m+1}} |h^{(m)}|_{L^{2}(\widetilde{I})}^{2}. $$
\end{lemma}
Applying Lemma 3-5 and the discussion in Section 3.1, the main theorem for the convergence of Poisson's equation can be obtained:
\begin{thm}
	Assuming that $u$ satisfies Eq. (\ref{eq:16}) and $u_{N}$ satisfies Eq. (\ref{eq:9}), if $u\in C^{m+2}(I)$ and $f\in C^{m}$, then there exists a constant $C:= C(L,\delta ,m, |u^{(m+2)}|_{L^{2}(I)} )$ such that
    \begin{equation}\label{eq:37}
    |u-u_{N}|_{L^{2}(I)} \leq CN^{-m}.
    \end{equation}
\end{thm}

As to the interpolation of convection-diffusion equations, the estimate of Lebesgue constant is still unsolved. More precisely, the following conjecture is of our concern.

\begin{conj}
	In addition to the assumtions given in Lemma 3, suppose that $u$ satisfies Eq. (\ref{eq:32}), and $u_{N}$ is defined by Eq. (\ref{eq:7}) such that
    \begin{equation}\label{eq:38}
    	-u_{xx}(x_{j})+k(x_{j})u_{x}(x_{j}) = -u_{N,xx}(x_{j})+k(x_{j})u_{N,x}(x_{j}), 
    \end{equation}
    for all $j=1,\ldots N$. Then, there exists a positive number $\sigma$ such that the associated Lebesgue constant $\Lambda_{N}$ satisfies
    \begin{equation}\label{eq:39}
    	\Lambda_{N} \leq CN^{\sigma},
    \end{equation}
    for all integers $N\geq 3$, where $C$ depends on $k$, $L$, and $\delta$.
\end{conj}

Once the Conjecture 1 is proved, the following conjecture can be obtained immediately.

\begin{conj}
	Inheriting the assumptions in Conjecture 1, there exists a constant $C:= C(L,\delta ,k,m, |u^{(m+2)}|_{L^{2}(I)})$ such that
    \begin{equation}\label{eq:40}
    	|u-u_{N}|_{L^{2}(I)} \leq CN^{-m+\sigma}.
    \end{equation}
\end{conj}

Although the estimate of the associated Lebesgue constant of convection-diffusion equations is still a mystery, our proposed method is applicable to most of cases. The discussion about this unsolved problem will be presented in Appendix. The way to avoid the difficulty is given in Appendix as well. Moreover, numerical test problems in Section 4 and 5 justify our proposed method.

\subsection{Generalization to higher dimensional problems and geometrically complex domains}

In order to estimate Poisson's equation and convection-diffusion equation of dimension $d\geq 2$, some classical results of PDE shall be used. Let $\Omega$ be a bounded domain in $R^{d}$, L is an second order elliptic operator such that
\begin{equation}\label{eq:41}
	Lu := \sum_{i,j=1}^{d} a^{ij}u_{x_{i}x_{i}}+\sum_{i=1}^{d}b^{i}u_{x_{i}}+cu = f.
\end{equation}
The following assumptions for $L$ are necessary: (A) There exists $\beta >0$ such that for all $x\in\Omega $, $\xi \in \mathbb{R^{d}}$, we have
\begin{equation}\label{eq:42}
	\sum_{i,j=1}^{d}a^{ij}(x)\xi_{i}\xi_{j} \geq \beta |\xi |^{2}.
\end{equation}
In addition, $a^{ij}(x)=a^{ji}(x)$ for all $i,j$, and $x$. (B) There exists a finite number $M$ such that
\begin{equation}\label{eq:43}
	\| a^{ij} \|_{C^{\alpha}(\Omega )} , \| b^{i} \|_{C^{\alpha}(\Omega )}, \| c \|_{C^{\alpha}(\Omega )} \leq M.
\end{equation}
Under the assumptions of $L$, the Schauder estimate for the cases of our concern can be given \cite{Jost2013}:

\begin{thm}
	Let $\Omega \subset \mathbb{R}^{d}$ be a bounded domain of class $C^{2,\alpha}$, $f\in C^{\alpha}(\overline{\Omega})$, $g\in C^{2,\alpha}(\overline{\Omega})$, $c(x)\leq 0$ and $u\in C^{2,\alpha}(\overline{\Omega})$, where $d\geq 2$ and $0<\alpha <1$, such that
    \begin{subequations}\label{eq:44}
    	\begin{align}
        	Lu = f &\ \ \ \ \text{in}\ \Omega , \\
            u=g    &\ \ \ \ \text{on}\ \partial\Omega .
        \end{align}
    \end{subequations}
    Then there exists a constant C that depends on $\Omega , \alpha , d, \beta$, and $M$ such that
    \begin{equation}\label{eq:45}
    	\| u \|_{C^{2,\alpha}(\Omega )} \leq C(\|f \|_{C^{\alpha}(\Omega )}+\|g \|_{C^{2,\alpha}(\Omega )}).
    \end{equation}
\end{thm}

Let R be the rectangular domain defined as in Section 2, and $f\in C^{m}_{0}(R)$. The argument in Lemma 2 and Lemma 5 are useful to $d$-dimensional cases, since the Fubini's theorem can be applied to derive
\begin{equation}\label{eq:46}
	\begin{aligned}
    	\int_{R}fw_{j_{1}\ldots j_{d}}d\mathbf{x} = \int_{I_{d}}\left[\cdots\int_{I_{2}}\left[\int_{I_{1}}f(x_{1},\ldots x_{d})w_{j_{1}}dx_{1} \right]dx_{2}\cdots \right]dx_{d}, 
    \end{aligned}
\end{equation}
where $I_{k}=(-\delta_{k},L_{k}+\delta_{k})$, $w_{j_{k}}=\sqrt{\frac{2}{L_{k}+2\delta_{k}}}\sin \left(\frac{j_{k}\pi(x_{k}+\delta)}{L_{k}+2\delta} \right)$ for $k=1,\ldots ,d$ and $j_{k}\in \mathbb{N}$. For simplicity, the uniform grid with $N$ grid points in each direction is considered. The d-dimensional trigonometric interpolation $f_{N}$ of $f$ with the basis defined by Eq. (\ref{eq:3}) can be expressed as:
\begin{equation}\label{eq:47}
	f_{N} = \sum_{j_{1},\ldots ,j_{d} = 1}^{N} f(x_{j_{1}},\ldots , x_{j_{d}})l_{j_{1}}(x_{1})\cdots l_{j_{d}}(x_{d}),
\end{equation}
since the variables of basis functions defined by Eq. (\ref{eq:3}) are separable. Now, all the above arguments for solving Poisson's equation can be integrated into a theorem:

\begin{thm}
	Let $R$ and $R_{0}$ be the domains defined in Section 2, $u\in C^{m+2,\alpha}(R_{0})$, $g\in C^{m+2,\alpha}(\overline{R_{0}})$, and $f\in C^{m,\alpha}(R_{0})$, $\alpha\in (0,1)$, where $u$, $g$, and $f$ satisfy Eq. (\ref{eq:1}), then there exists a constant $C:= C(R,m,| u |_{H^{m+2}(R_{0})})$ such that
    \begin{equation}\label{eq:48}
    	|u-u_{N}|_{L^{2}} \leq CN^{-m},
    \end{equation}
where $| u |_{H^{m+2}(R_{0})} = \sup_{|\gamma| = m+2}|D^\gamma u|_{L^{2}(R_{0})} $.
\end{thm}
Moreover, once Conjecture 1 holds, we have the following:
\begin{conj}
	Let $R$ and $R_{0}$ be the domains defined in Section 2, $u\in C^{m+2,\alpha}(R_{0})$, $k\in [C^{m,\alpha}(R_{0})]^{d}$, $g\in C^{m+2,\alpha}(\overline{R_{0}})$, and $f\in C^{m,\alpha}(R_{0})$, $\alpha\in (0,1)$, where $u,k,g$, and $f$ satisfy
    \begin{subequations}\label{eq:49}
    	\begin{align}
        	-\Delta u + k\cdot \nabla u &= f  \ \ \ \  \text{in}\ R \\
            u &= g \ \ \ \ \text{on} \ \partial R.
        \end{align}
    \end{subequations}
    Then there exists a constant $C:= C(R,k,m,| u |_{H^{m+2}(R_{0})})$ such that
    \begin{equation}\label{eq:50}
    	|u-u_{N}|_{L^{2}} \leq CN^{-m+\sigma}
    \end{equation}
    for some fixed $\sigma >0$
\end{conj}

To deal with geometrically complex domains, the strategy of moving grid points is introduced in Section 2.2. We can find the fact that \textit{every moved grid point will not be apart from its original position over half of the grid spacing.} That is, let $x_{j}$ be a grid point and $x_{j}^{'}$ its new position after the grid movement induced by fitting the boundary of the physical domain, we must have
\begin{equation}\label{eq:51}
	|x_{j}^{'}-x_{j}| \leq \frac{1}{2}h,
\end{equation}
where $h = \max_{j=1,\ldots d} L_{j}/N$. Recalling the arguments from Eq. (\ref{eq:19}) to Eq. (\ref{eq:30}), Theorem 1 still holds if the new grid points are chosen as Eq. (\ref{eq:51}). Therefore, the convergence for the approximation of Poisson's equation in geometrically complex domain can be ensured. As to the convection-diffusion equation, the convergence can be achieved once Conjecture 1 is proved. More details of dealing with the convection-diffusion equation in geometrically complex domain are discussed in Appendix.

\section{Numerical results of two-dimensional problems with rectangular domain}
\subsection{Poisson's equation}
\begin{figure}\label{simple_poisson}
	\centering
	\includegraphics[width=\textwidth]{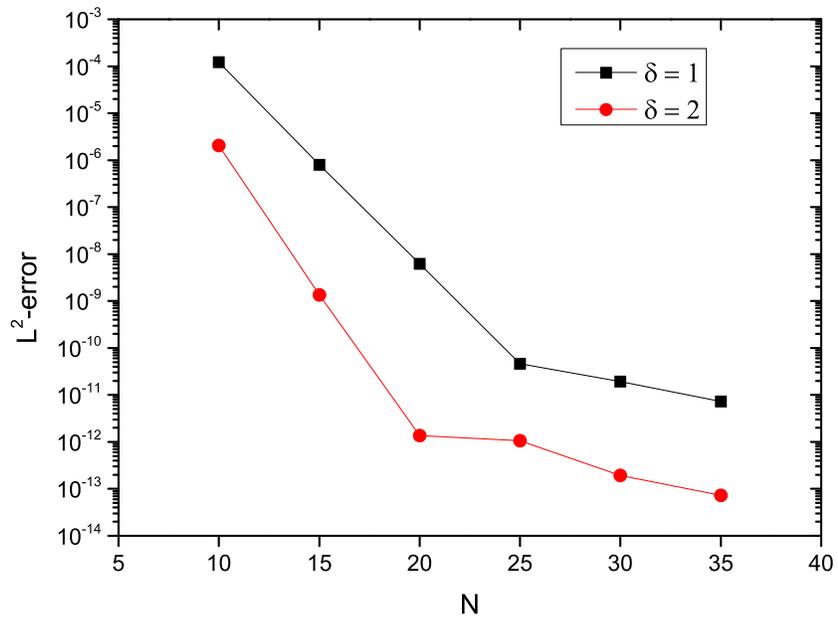}
    \caption{The results show the influence of the stretching number $\delta$. In this case, we consider $\delta := \delta_{1}=\delta_{2}$ defined in Section 2.}
\end{figure}

We firstly verify our proposed method applied to solve Poisson's equation. Let $\Omega = (0,2)\times (0,2)$, and $u_{e}=x^{2}y^{3}$. The Poisson's equation in Eq. (\ref{eq:1}) with $R_{0}=\Omega$ is considered. Furthermore, $f$ and $g$ are chosen properly so that $u_{e}$ is the exact solution. The results is given in Fig. \ref{simple_poisson}. We find that stretching the basis functions can indeed affect the convergence. 

\subsection{Convection-diffusion equation}
The test problem for convection-diffusion equation is given by Chiu et al. \cite{CHIU20093640}, which is suitable for the verification. Consider the convection-diffusion in the domain $\Omega=(0,1)\times(0,1)$ such that
\begin{equation}\label{eq:52}
	u\frac{\partial\phi}{\partial x}+v\frac{\partial\phi}{\partial y}=\frac{1}{Re}\left( \frac{\partial^{2}\phi}{\partial x^{2}}+\frac{\partial^{2}\phi}{\partial y^{2}}\right) +S.
\end{equation}
In the above, $Re$ and $S$ are the Reynolds number and the source term, respectively. Moreover,
\begin{subequations}\label{eq:53}
	\begin{align}
    	u = \frac{-2(1+y)}{(1+x)^{2}+(1+y)^{2}} \\
         v = \frac{2(1+x)}{(1+x)^{2}+(1+y)^{2}}.
    \end{align}
\end{subequations}
In this test problem, we set the exact solution $\phi$ to be identical to $u$ given in Eq. (\ref{eq:53}a), provided that $S=-\frac{\partial p}{\partial x}$, where
\begin{equation}\label{eq:54}
	p = -\frac{2}{(1+x)^{2}+(1+y)^{2}}.
\end{equation}
Now, we are going to show the good spatial rates of convergence by using our proposed scheme in the above test problem. The $L^{2}$-error norms computed at $10\times 10$, $15\times 15$, $20\times 20$, $25\times 25$, and $30\times 30$ uniform grids are performed. In view of Fig. \ref{simple_re}, our proposed method exhibits rapid convergence for solving the convection-diffusion equation as well.
\begin{figure}
	\centering
    \subfigure[$Re=10$]{
    \label{re10}
    \includegraphics[width=0.8\linewidth]{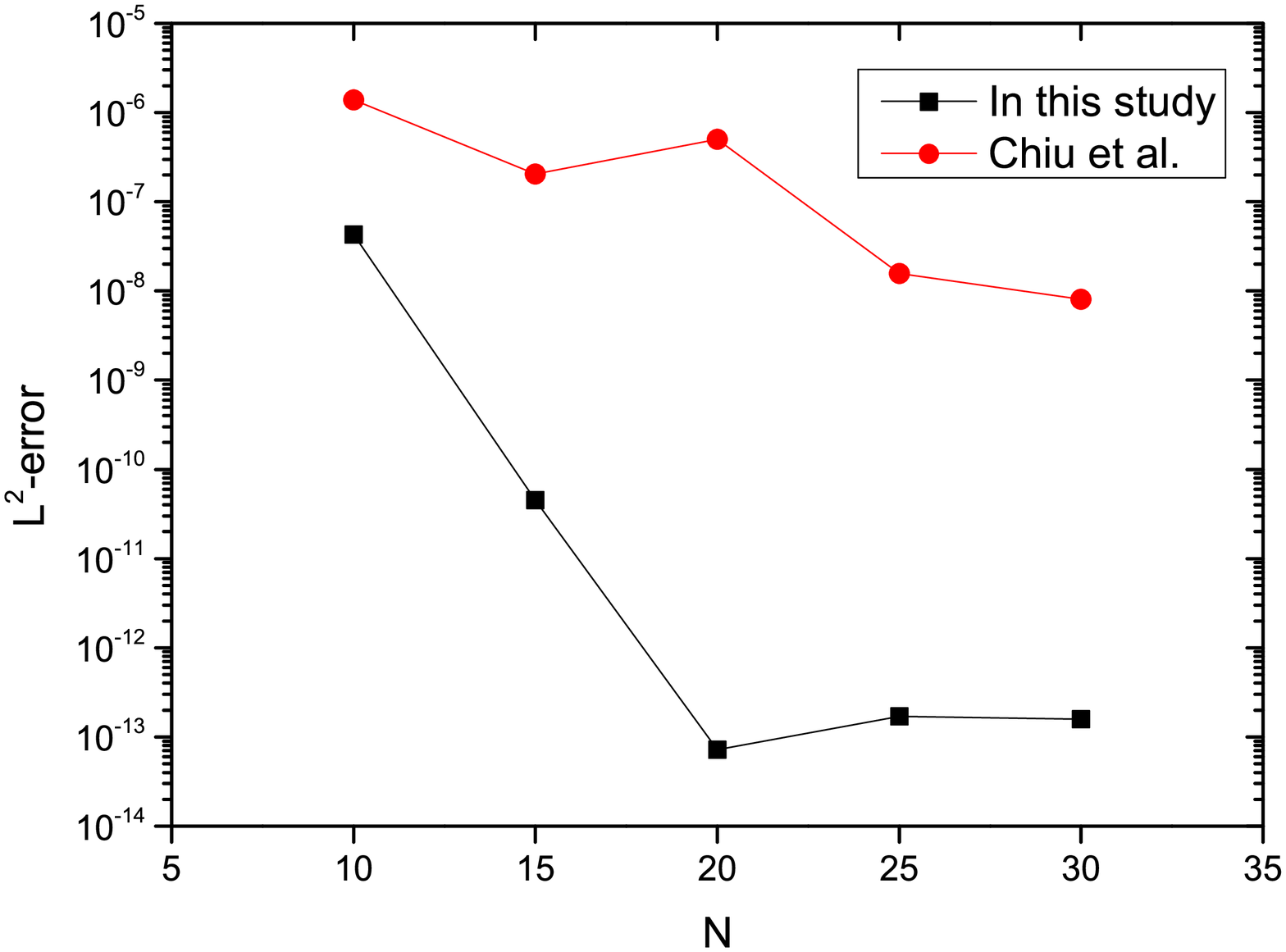}
    }
    \subfigure[$Re=100$]{
    \label{re100}
    \includegraphics[width=0.8\linewidth]{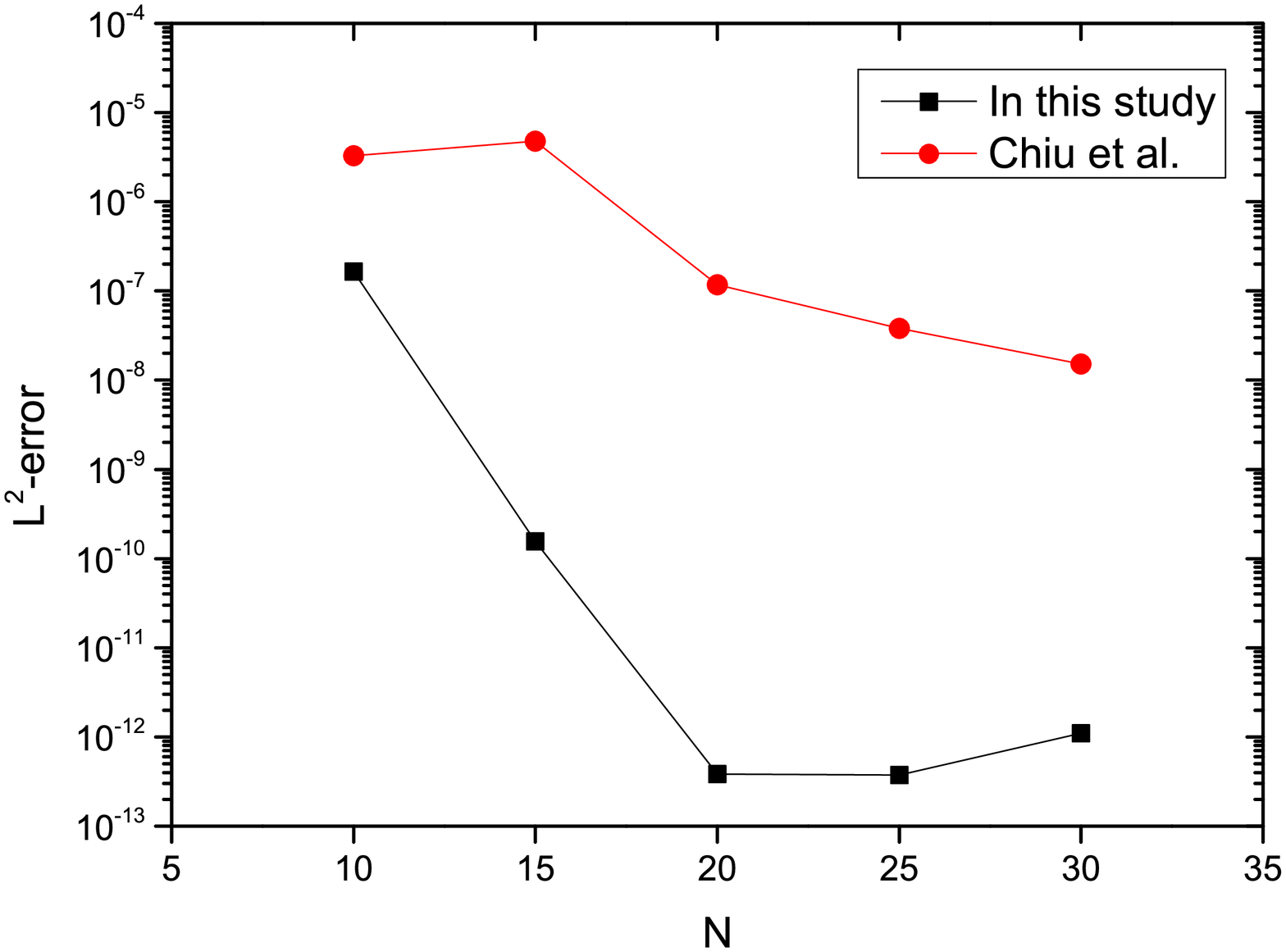}
    }
    \captcont{Comparison of the $L^{2}$-error norms. Reynolds numbers $Re=10,100,10000,100000$ are considered in this study. The stretching number $\delta$ is chosen to be $2$ in each case.}
    \label{simple_re}
\end{figure}

\begin{figure}
\centering
    \subfigure[$Re=10000$]{
    \label{re10000}
    \includegraphics[width=0.8\linewidth]{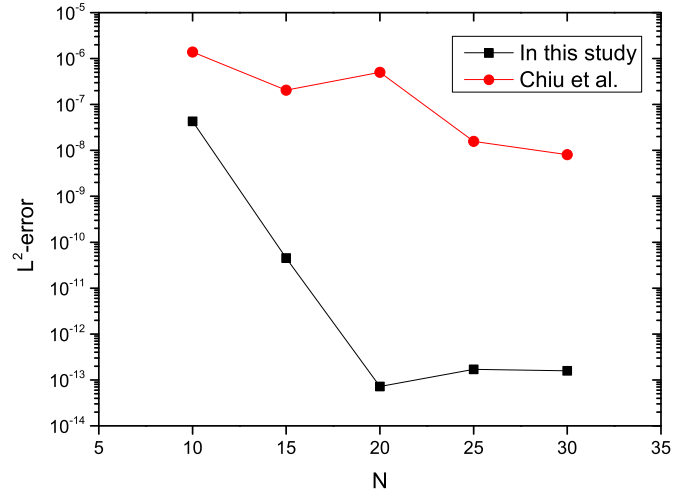}
    }
    \subfigure[$Re=100000$]{
    \label{re100000}
    \includegraphics[width=0.8\linewidth]{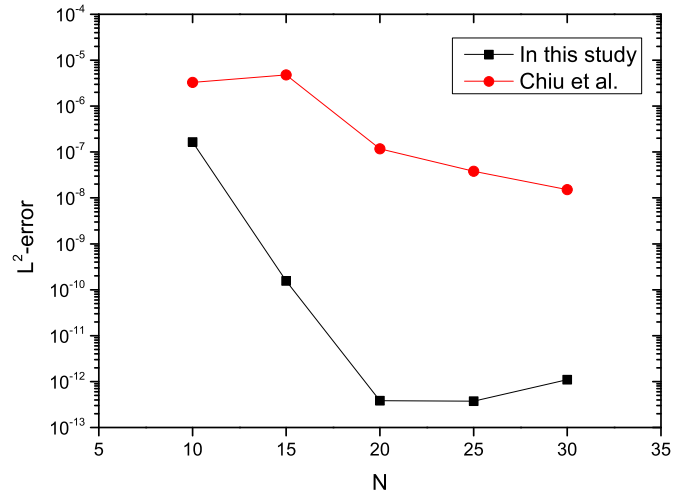}
    }
    \caption*{Comparison of the $L^{2}$-error norms. Reynolds numbers $Re=10,100,10000,100000$ are considered in this study. The stretching number $\delta$ is chosen to be $2$ in each case. (cont.)}
    \label{simple_re}
\end{figure}

\begin{figure}
	\centering
	\includegraphics[width=\textwidth]{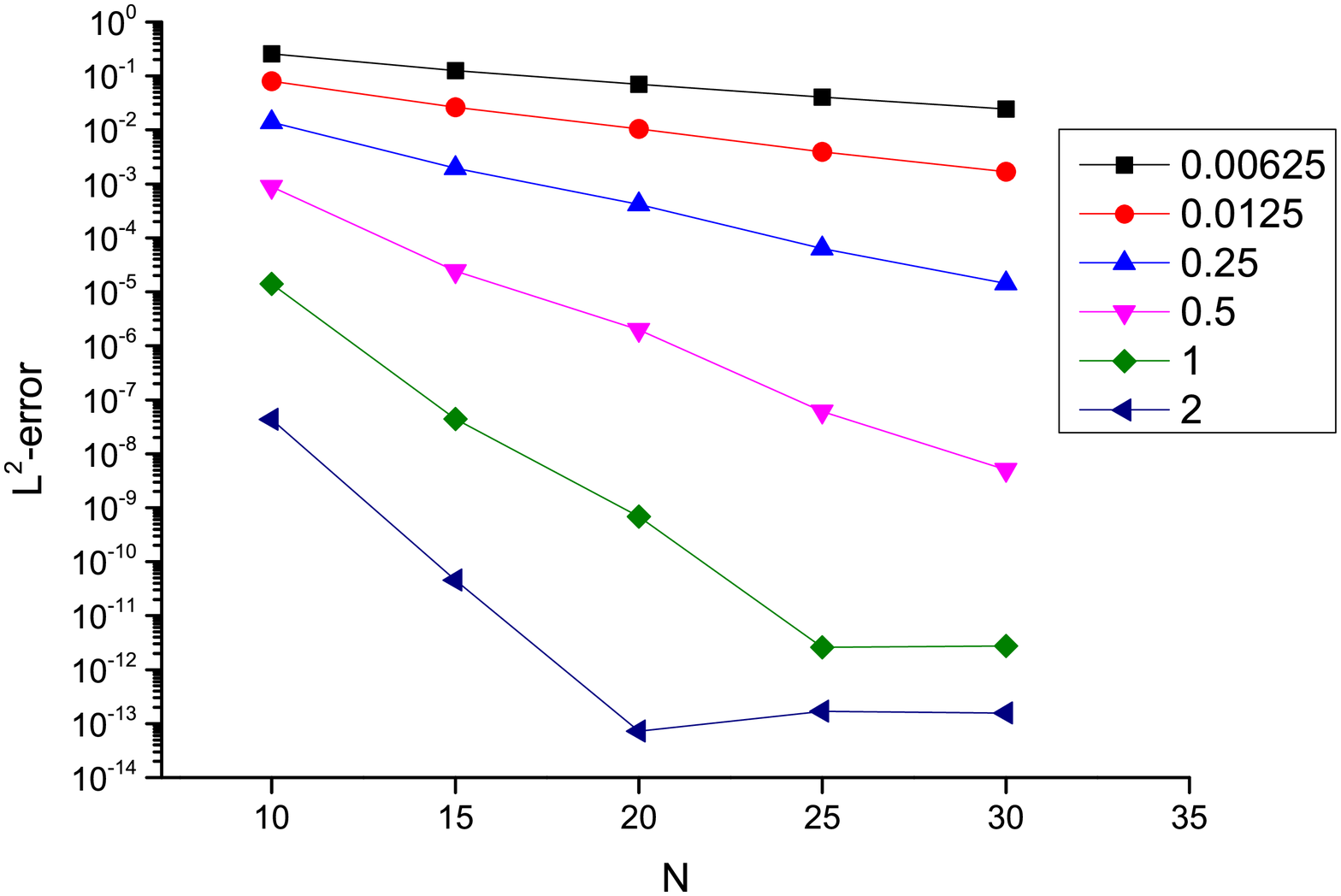}
    \caption{Comparison of the $L^{2}$-error norms at fixed $Re=10$ and at different stretching numbers, which are $\delta = 0.0625,0.125,0.25,0.5,1,2$.}
    \label{simple_cd_compared}
\end{figure}

Similar to the case of solving the Poisson's equation, the choice of the stretching number $\delta$ is a crucial issue here. The effect of the stretching number $\delta$ can be clearly seen in Fig. \ref{simple_cd_compared}.

\section{Numerical results of two-dimensional problems with geometrically complex domain}
In this section, the predicted errors will be cast in $L^{\infty}$ norm. We need the following results derived by Morrey's inequality \cite{Evans2010}
\begin{lemma}
	Let $\Omega$ be a bounded domain in $R^{d}$, $d\in \mathbb{N}$, $u\in H_{0}^{1}(\Omega)$. Then there exists a constant $C$ depending on $d$ and $\Omega$ such that
    \begin{equation*}
    	\|u \|_{C^{0,1-\frac{d}{2}}(\Omega )} \leq C \|u \|_{H^{1}(\Omega )},
    \end{equation*}
    where $\|\cdot \|_{C^{0,\alpha}(\Omega )}$ is the H$\ddot{o}$lder norm for $0<\alpha <1$.
\end{lemma}
This ensures the convergence of $L^{\infty}$-error if the solution of PDEs possesses sufficient smoothness.
\subsection{Poisson's equation}
\begin{figure}
	\centering
    \includegraphics[width=\textwidth]{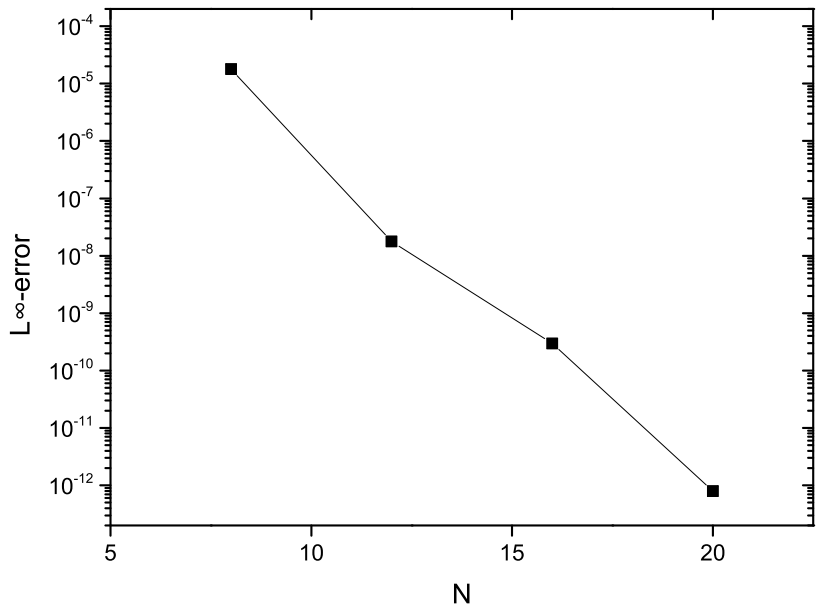}
    \caption{The predicted $L^{\infty}$-error norms against $N$ for Eq. (\ref{eq:55}).}
    \label{ibse}
\end{figure}

To show that the proposed method is applicable to geometrically complex cases, the Poisson's equation considered by Stein et al. \cite{STEIN2016252} in the domain $\Omega = B_{2}((\pi , \pi ))$ is investigated:
\begin{subequations}\label{eq:55}
	\begin{align}
    	-\Delta u = 4,\ \ \ & \text{in } \Omega , \\
        u=0,\ \ \ & \text{on } \partial \Omega .
    \end{align}
\end{subequations}
The exact solution $u_{e}$ is given by
\begin{equation}\label{eq:56}
	u_{e} = 4-(x-\pi )^{2} - (y-\pi )^{2}
\end{equation}
The numerical results show that the high accuracy nature of the spectral method is still retained using our proposed method, see Fig. \ref{ibse}. The approximation solution of Eq. (\ref{eq:55}) generated by IBSE method cannot reach such convergence since the boundary smooth extension in IBSE method restricts the regularity of the original PDE. However, such a limitation won't happen in our proposed method.

\begin{figure}
	\centering
    \subfigure[The fish domain.]{
    \label{fish_domain}
    \includegraphics[width=0.6\linewidth]{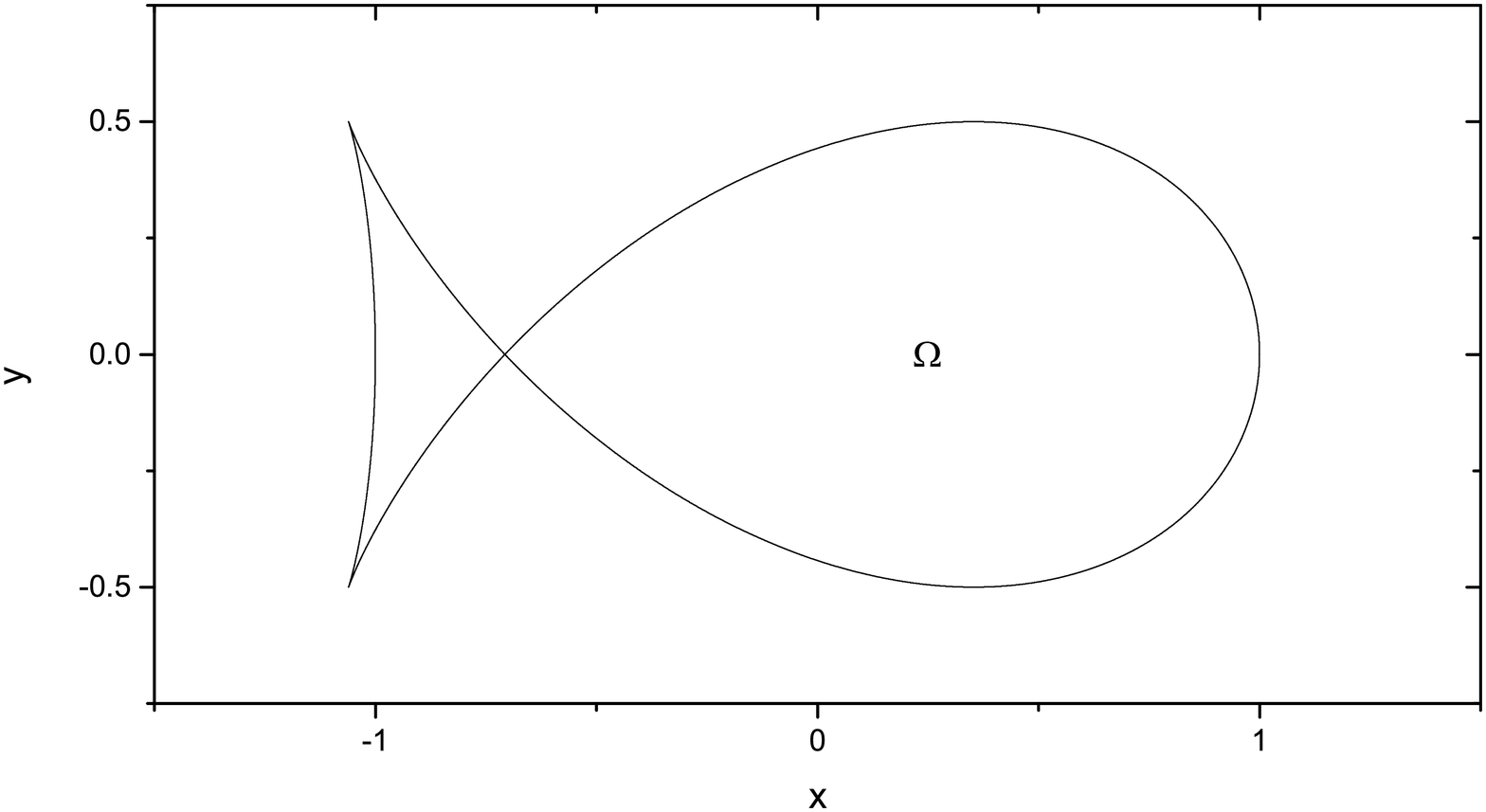}
    }
    \subfigure[Numerical results.]{
    \label{fish_result}
    \includegraphics[width=\linewidth]{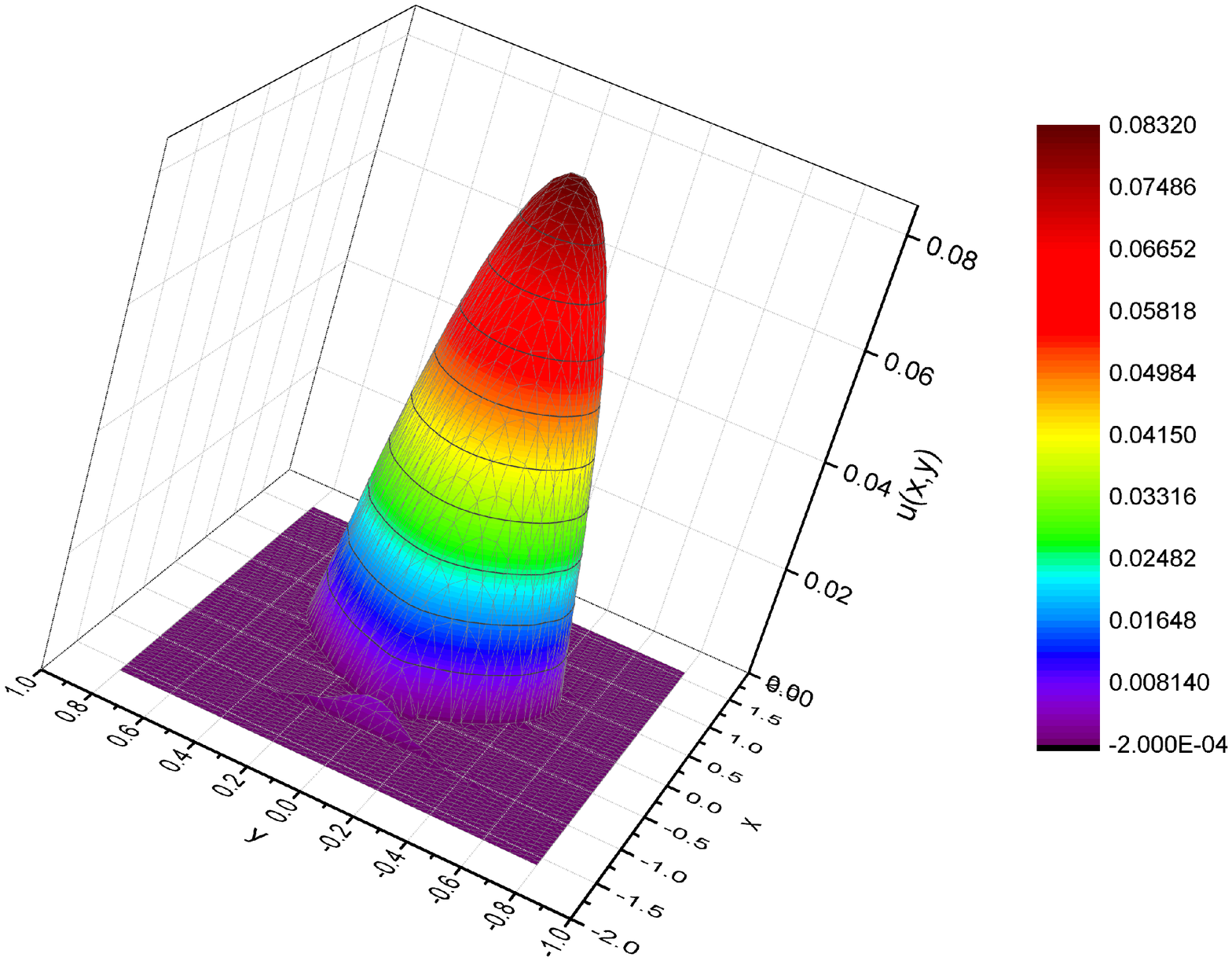}
    }
    \caption{The domain $\Omega$ is defined by a closed curve, which can be expressed as $(2x^{2}+y^{2})^{2}-2\sqrt{2}x(2x^{2}-3y^{2})+2(y^{2}-x^{2})=0$. }
    \label{fish}
\end{figure}

\begin{figure}
	\centering
    \subfigure[The ''ninja dart'' domain.]{
    \label{ninja_domain}
    \includegraphics[width=0.4\linewidth]{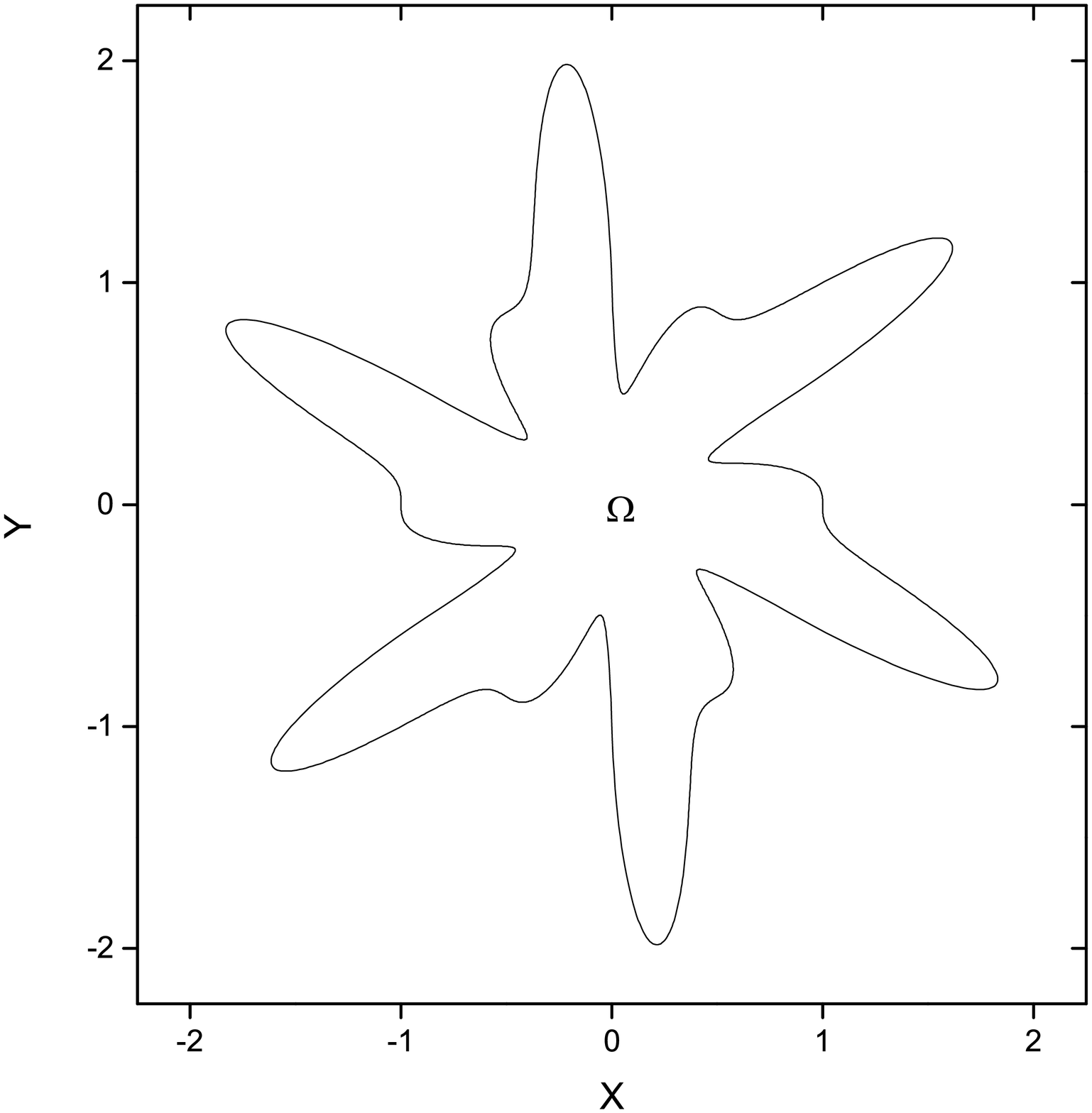}
    }
    \subfigure[Numerical results.]{
    \label{ninja_result}
    \includegraphics[width=\linewidth]{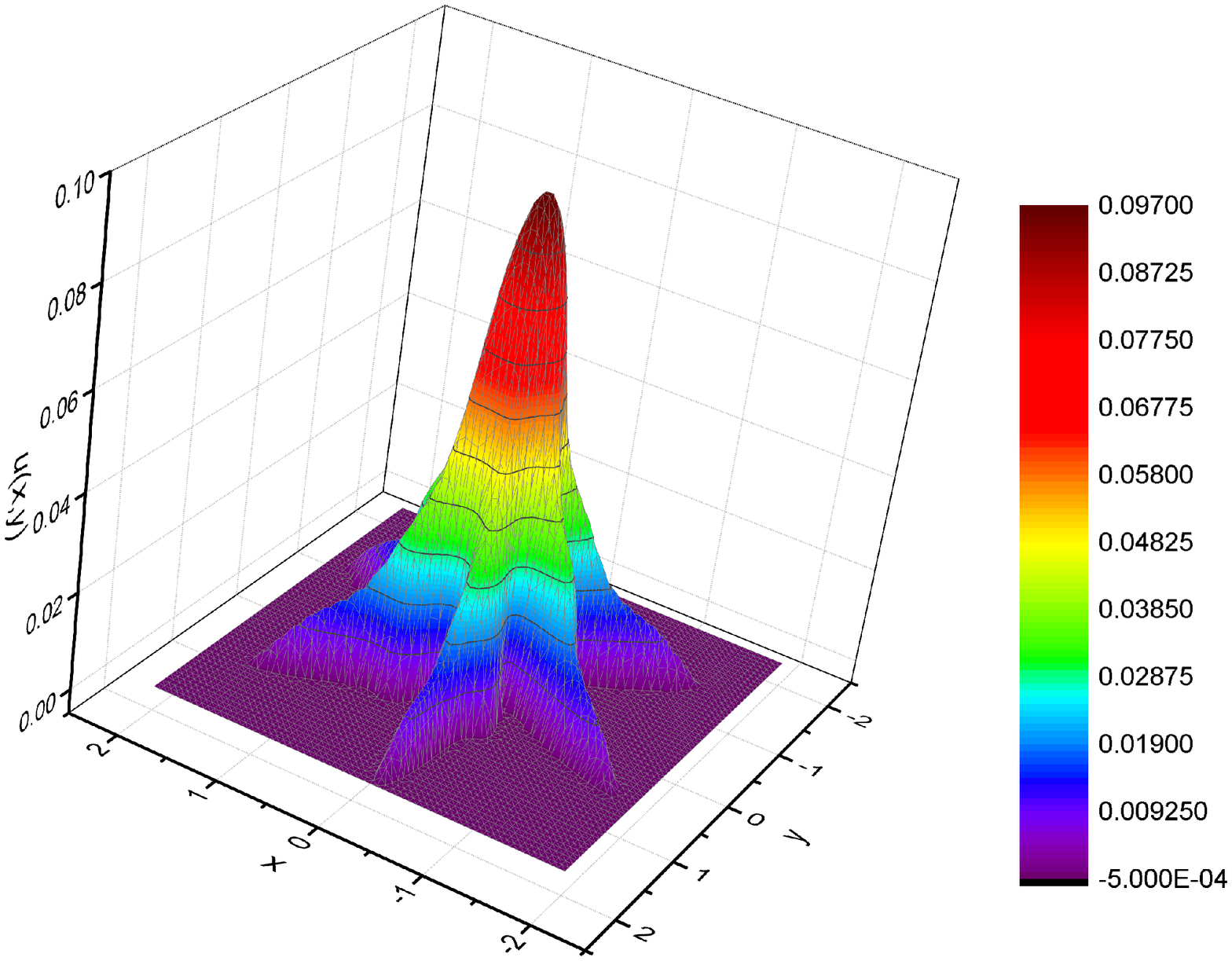}
    }
    \caption{The domain $\Omega$ is defined by a closed curve, which can be expressed as, in polar coordinate, $r = |\cos(3t)|^{\sin(6t)}$, for $t\in [0,2\pi)$.}
    \label{ninja}
\end{figure}

We will not be satisfied if merely the numerical results of a two-dimensional ball are shown. Therefore, another two numerical examples with different geometries are given. In both cases, the following governing equation is considered and the domains $\Omega$ are defined in Fig. \ref{fish} and \ref{ninja}. 
\begin{subequations}\label{eq:add1}
	\begin{align}
		-\Delta u = 1, \ \ \ & \text{in } \Omega , \\
		u = 0, \ \ \ & \text{on } \partial \Omega.
	\end{align}
	
\end{subequations}

\subsection{Convection-diffusion equation}
\begin{figure}
	\centering
    \subfigure[For $a=0.9$ and $b=0.6$]{
    \label{elliptic_result}
    \includegraphics[width=0.8\linewidth]{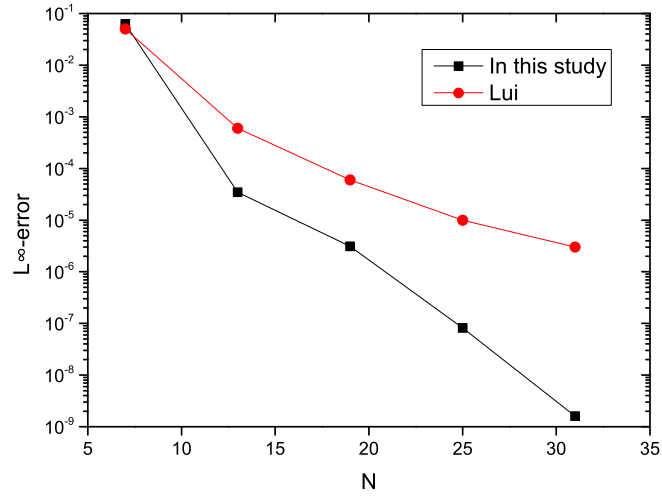}
    }
    \subfigure[For $a=0.9$ and $b=0.9$]{
    \label{circle_result}
    \includegraphics[width=0.8\linewidth]{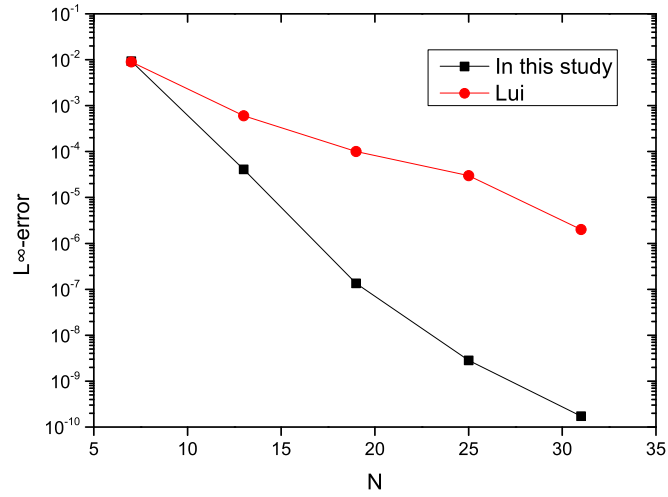}
    }
    \caption{The predicted onvergence rate ($L^{\infty}$-error norms) for the test problem in Section 5.2.}
    \label{elliptic_compared}
\end{figure}

Two test problems of Lui \cite{LUI2009541} used for the sake of verification are investigated here. For the first case, the elliptic domain is considered. Let $\Omega = \lbrace (x,y)\ |\ \Gamma (x,y)<0 \rbrace $, where $\Gamma (x,y)=\frac{x^{2}}{a^{2}}+\frac{y^{2}}{b^{2}}-1 $. The convection-diffusion equation under investigation is given by
\begin{subequations}\label{eq:57}
	\begin{align}
    	-u_{xx}-u_{yy}-u_{x}-u_{y}=f,& \ \ \text{in} \ \Omega, \\
        u = 0, & \ \ \text{on} \ \partial \Omega.
    \end{align}
\end{subequations}
In the above, $f$ is properly chosen so that the exact solution $u_{e}$ takes the form of
\begin{equation}\label{eq:58}
	u_{e} = \sin \left( \frac{\pi}{2}\Gamma (x,y) \right).
\end{equation}
The numerical results are shown in Fig. \ref{elliptic_compared}.
\begin{figure}
	\centering
    \subfigure[$a=0.9$ and $b=0.6$]{
    \label{elliptic}
    \includegraphics[width=0.6\linewidth]{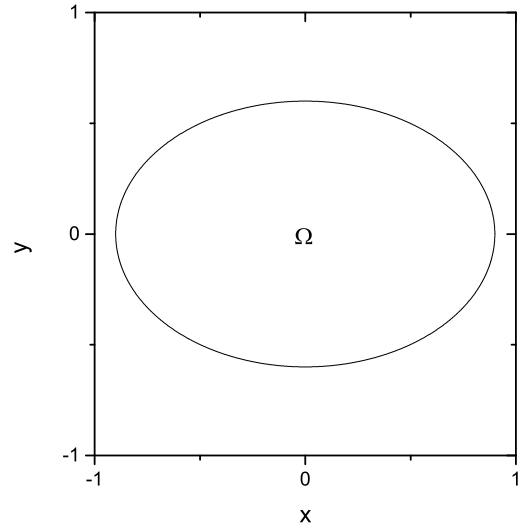}
    }
    \subfigure[$a=0.9$ and $b=0.9$]{
    \label{circle}
    \includegraphics[width=0.6\linewidth]{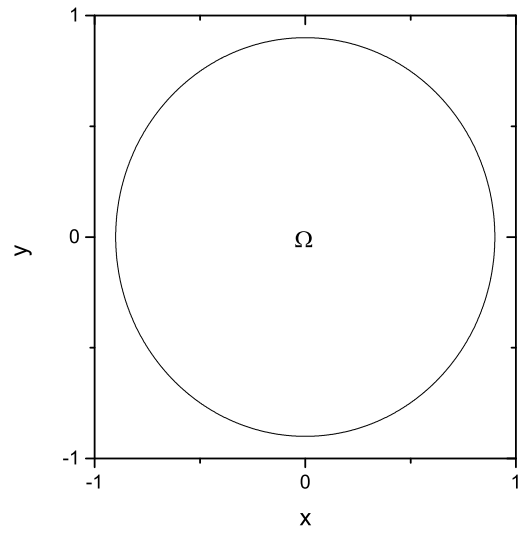}
    }
    \caption{The physical domain $\Omega$ considered for the test problem considered in Section 5.2}
    \label{elliptic_domain}
\end{figure}

\begin{figure}
	\centering
    \subfigure[The star-shaped domain]{
    \label{star_domain}
    \includegraphics[width=0.6\linewidth]{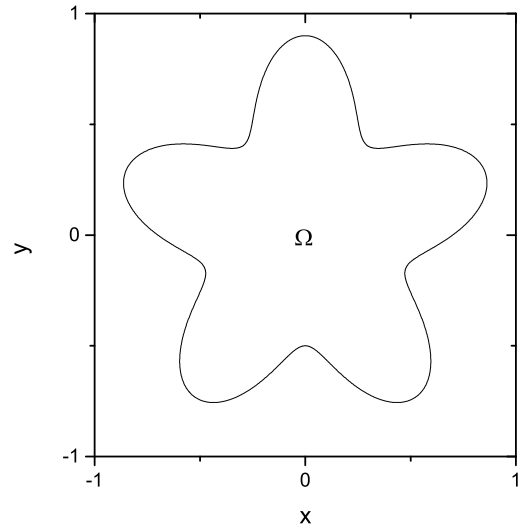}
    }
    \subfigure[The $L^{\infty}$-error norms]{
    \label{star_result}
    \includegraphics[width=0.8\linewidth]{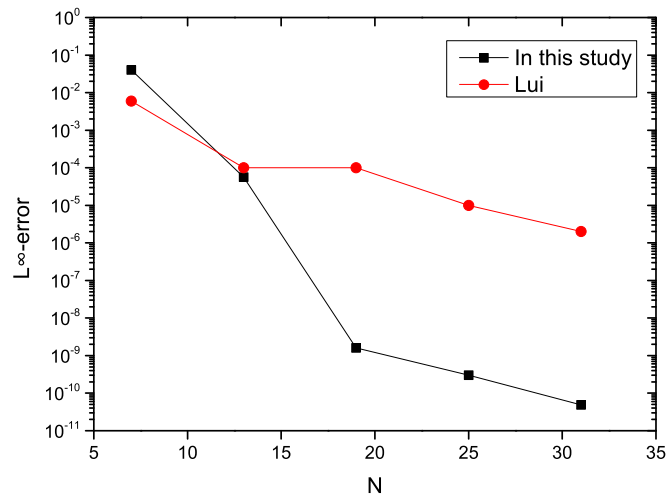}
    }
    \caption{Numerical result for the test problem in Section 5.2}
    \label{star_compared}
\end{figure}

The second test problem is the convection-diffusion equation in a star-shaped domain. Let $\Omega$ be the domain such that its boundary is defined by the polar equation $r=0.7+0.2\sin 5\theta $. The equation to be solved is the same as Eq. (\ref{eq:57}a), except that the Dirichlet boundary condition and $f$ are chosen properly so that the exact solution is $u_{e}= \sin (x^{2}+y^{2})$. The numerical result is plotted in Fig. \ref{star_compared}.

\section{Conclusion and future works}
In this study, we have developed a brand new method for solving PDEs in the framework of spectral methods. The method is shown to be as accurate as the previously proposed spectral method for geometrically simple problem. Of greater significance, the method can reach a much higher accuracy than others in geometrically complex problems. In our proposed method, the chosen basis functions can be simply regarded as the extension of sine functions. We extend the domain within which basis functions are defined so that the solutions of PDEs can be calculated with an extremely high accuracy. The reason for rendering high accuracy of our proposed method is given. In comparison with the spectral collocation methods proposed in the past, the accuracy of our scheme is not sensitive to the choice of collocation points. The error caused by collocation can be eliminated by the extension of basis functions. Given the advantage of the high degree of freedom for collocating points, our proposed method enjoys high accuracy as well in the simulation of geometrically complex problems.

The crucial role of the convergence rate in our proposed method is Lebesgue constant. We have shown that for the Poisson's equation in a rectangular domain, the error can be well reduced by stretching the domain for basis functions. The Lebesgue constant for solving the Poisson's equation can be bounded by an independent constant once a proper stretching of basis functions is conducted. High accuracy of solving Poisson's equation in a geometrically complex domain can be explained by the upper bound of the Lebesgue constant as well. When solving geometrically complex problems, the moved grid points won't be distant from the original position too far, our theory for the upper bound of the Lebesgue constant is still valid. Therefore, in all cases, the rapid convergence rate in solving the Poisson's equation is confirmed. As to the simulation of convection-diffusion equation, the best approximation in view of the stretched basis functions is shown in Section 3. The only problem of the convergence rate in solving the convection-diffusion equation is, similar to the Poisson's equation, the determination of Lebesgue constant.  In Appendix A, we show that the trigonometric Lagrange polynomial can be expressed in terms of the determinants of the matrices containing the values of basis functions on each collocation point. However, the formula of the Lagrange polynomials are too complicated to get the analytical upper bound of the Lebesgue constant associated to the convection-diffusion equation. This difficulty can be partially resolved by calculating the Lebesgue constant at each kind of fixed grid, shown in Fig. \ref{fig:a1}. The error can be eliminated by avoiding some improper grids and the stretched basis functions.

Several two-dimensional numerical results are investigated in this study. For the Poisson's equation and the convection-diffusion equation in the rectangular domain, numerical experiments show that our method  retains the exponentially decay nature as the traditional spectral method does. This shows, at least, our method won't be degraded for solving PDEs with simple geometry. As to the geometrically complex problem, our proposed method is proved to be much better than other methods. On the test of the convection-diffusion equation, our method is at least 3 to 4 orders better than Lui \cite{LUI2009541} for each test problem. For the Poisson's equation, exponentially decay convergence nature is still retained. We need only 20 points in each direction to achieve the $L^{\infty}$-error norm of $10^{-12}$. Furthermore, geometrically complicated problems can be handled using the proposed method.    

The limitation of this work is that the high-accuracy cannot be expected if the regularity of the original PDE is low. We shall extend the method in this study to different kinds of boundary condition. Whether our proposed method is suitable for nonlinear PDEs and time-dependent problems is worth studying in the future.

\appendix
\section{Estimates of the Lebesgue constant for convection-diffusion equation}

We firstly consider the convection-diffusion equation Eq. (\ref{eq:32}) and define $\psi_{j}(x)$ such that
\begin{equation}\label{eq:a1}
	\psi_{j}(x) = \sqrt{\frac{2}{L+2\delta}}\sin \left(\frac{j\pi(x+\delta)}{L+2\delta} +\theta_{j}(x) \right),
\end{equation}
where
\begin{equation}\label{eq:a2}
	\theta_{j}(x) =
    \begin{cases}
    \sin^{-1}\left(\frac{1}{\sqrt{1+\left(\frac{j\pi}{k(x)} \right)^{2}}}  \right) & \text{if}\ k(x) \neq 0 \\
    0 & \text{if}\ k(x) = 0.
    \end{cases}
\end{equation}
Let $f_{N}$ be the interpolation of $f$ defined in Lemma 4, $f_{N}$ can be expressed as
\begin{equation}\label{eq:a3}
	f_{N}(x)=\sum_{j=1}^{N}f(x_{j})l_{j}(x).
\end{equation}
In the above, $l_{j}(x)$ can be expressed in terms of the basis functions $\lbrace \psi_{j}(x)\rbrace_{j=1}^{N} $ defined in Eq. (\ref{eq:a1}) for the one-dimensional case
\begin{equation}\label{eq:a4}
	l_{j}(x) = \frac{1}{\det \Psi}\sum_{i=1}^{N}(-1)^{i+j}\psi_{i}(x)\det \Psi_{ji},
\end{equation}
where $\Psi:=(\psi_{i}(x_{j}))_{i,j=1}^{N}$, $\Psi_{ji}$ is the $ji$-cofactor of $\Psi$, and $x_{j}$'s are the grid points. It is easy to find that
\begin{equation}\label{eq:a5}
		l_{j}(x_{k})= 
        \begin{cases}
        1 & \text{if}\ k = j, \\
    	0 & \text{if}\ k \neq j.
        \end{cases}
\end{equation}
Therefore, Eq. (\ref{eq:a3}) is indeed a Lagrange interpolation of $f$. The estimates of Lebesgue constant for such interpolation can be estimated by finding the upper bound of $|l_{j}(x)|$ for $j=1,\ldots ,N$. In Section 3, we have mentioned the difficulty of the estimate of Lebesgue constants and regard it as a conjecture. Nonetheless, the difficulty can be partially resolved by doing some numerical exercises in advance.

\begin{figure}
	\includegraphics[width=\textwidth]{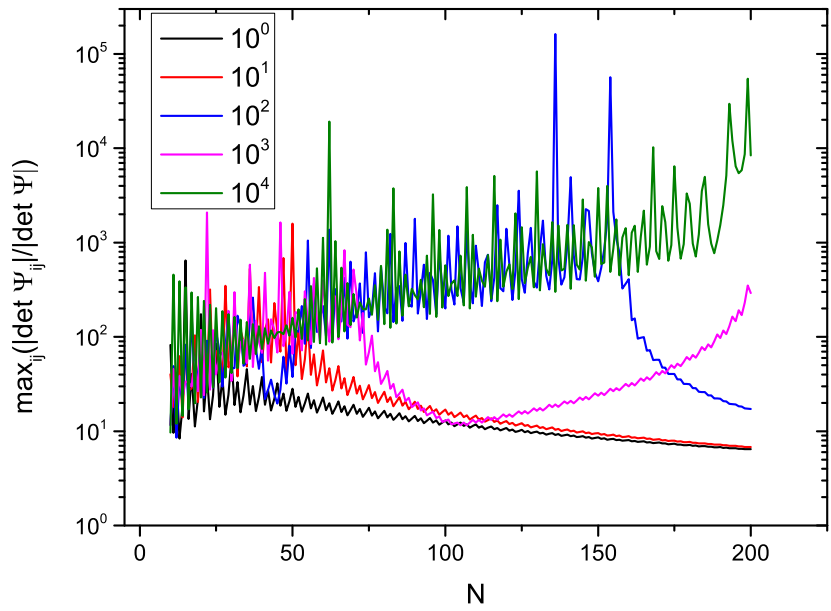}
    \caption{Boundedness of the Lebesgue constant against different constant convection  $k=1,10,100,1000,10000$. Here $\delta=2.01$ is considered.}\label{fig:a1}
\end{figure}

In view of Eq. (\ref{eq:a4}), we have
\begin{equation}\label{eq:a6}
	|l_{j}(X)| \leq \frac{1}{|\det \Psi|}\sum_{i=1}^{N}|\det \Psi_{ji}|,
\end{equation}
since $|\psi_{i}(x)|\leq 1$. The upper bound of the Lebesgue constant $\Lambda_{N}:=\sup_{x}\sum_{j=1}^{N}|l_{j}(x)|$ for the fixed $N$ can be investigated by the ratio between the absolute values of $\Psi$ and $\Psi_{ij}$ for $i,j=1,\ldots N$. We may find a proper $\delta$ such that the Lebesgue constant $\Lambda_{N}$ won't be a large number. Therefore, the problem of solving the convection-diffusion equation can be treated by plotting a graph, for example, Fig. \ref{fig:a1}. This helps us to determine a proper $\delta$ so that the boundedness of the Lebesgue constant can be ensured.

\bibliographystyle{model1-num-names}







\end{document}